\documentclass{amsart}

\usepackage{dsfont}
\usepackage{mathtools}
\usepackage{enumitem}
\usepackage{amsmath}
\usepackage{amssymb}

\numberwithin{equation}{section}

\newtheorem{theorem}{Theorem}[section]
\newtheorem{corollary}[theorem]{Corollary}
\newtheorem{lemma}[theorem]{Lemma}
\newtheorem{proposition}[theorem]{Proposition}

\theoremstyle{definition}
\newtheorem{definition}[theorem]{Definition}
\newtheorem{example}[theorem]{Example}
\newtheorem{remark}[theorem]{Remark}

\newcommand{\eps}{\varepsilon}	 
\newcommand{\defeq}{\vcentcolon=}
\newcommand{\eqdef}{=\vcentcolon}
\newcommand{\inte}{\textup{int}}
\newcommand{\om}{\omega}
\newcommand{\ee}{\mathrm{e}}
\newcommand{\mdim}{\delta}
\newcommand{\leb}{\lambda}
\newcommand{\bij}{\pi }
\newcommand{\bd}{\rho}
\newcommand{\bdneu}{p}

\newcommand{\card}{\#}
\newcommand{\eigenf}{h}
\newcommand{\eigenv}{\gamma}
\newcommand{\diam}{\textup{diam}}

\newcommand{\omneu}{u}

\newcommand{\main}{image}
\newcommand{\ze}{\zeta}
\newcommand{\aaa}{a}
\newcommand{\bb}{b\,}
\newcommand{\bbb}{b}

\newcommand{\entro}{H}
\newcommand{\sset}{K}

\renewcommand{\theenumi}{(\roman{enumi})}

\DeclareMathOperator*{\wlim}{w-lim\,}

\title[Minkowski Content and Fractal Euler Characteristic for cGDS]{Minkowski Content and Fractal Euler Characteristic for Conformal Graph Directed Systems}

\author[M.~Kesseb\"ohmer]{Marc Kesseb\"ohmer}
\address{M.~Kesseb\"ohmer, Universit\"at Bremen, Bibliothekstra{\ss}e 1, 28359 Bremen, Germany}
\email{mhk@math.uni-bremen.de}

\author[S.~Kombrink]{Sabrina Kombrink}
\address{S.~Kombrink, Universit\"at Bremen, Bibliothekstra{\ss}e 1, 28359 Bremen, Germany}
\curraddr{Universit\"at zu L\"ubeck, Ratzeburger Allee 160, 23562 L\"ubeck, Germany}
\email{kombrink@math.uni-luebeck.de}

\thanks{The second author was supported by grant 03/113/08 of the Zentrale Forschungsf\"orderung, Universit\"at Bremen.}


\begin{document}

\begin{abstract}
	We study the  (local) Minkowski content and the (local) fractal Euler characteristic of  limit sets $F\subset\mathbb R$ of conformal graph directed systems (cGDS) $\Phi$. For the local quantities we prove that the logarithmic Ces\`aro averages always exist and are constant multiples of the $\mdim$-conformal measure. If $\Phi$ is non-lattice, then also the non-average local quantities exist and coincide with their respective average versions. When the conformal contractions of $\Phi$ are analytic, the local versions exist if and only if $\Phi$ is non-lattice.
	For the non-local quantities the above results in particular imply that limit sets of Fuchsian groups of Schottky type are Minkowski measurable, proving a conjecture of Lapidus from 1993.
  Further, when the contractions of the cGDS are similarities, we obtain that the Minkowski content and the fractal Euler characteristic of $F$ exist if and only if $\Phi$ is non-lattice, generalising earlier results by Falconer, Gatzouras, Lapidus and van Frankenhuijsen for non-degenerate self-similar subsets of $\mathbb R$ that satisfy the open set condition.
\end{abstract}

\maketitle

\section{Introduction}\label{sec:intro}

We examine the (local) Minkowski content and the (local) fractal Euler characteristic of limit sets of finite conformal graph directed systems (cGDS) that are embedded in $\mathbb R$, as introduced e.\,g.\ in \cite{Urbanski_Buch}. 
The class of cGDS gives rise to a rich collection of fractal sets including self-conformal sets, limit sets of Fuchsian groups of Schottky type and limit sets of Markov interval maps. 
So far, the (local) Minkowski content and the (local) fractal Euler characteristic have been investigated only for restrictive subclasses of limit sets of cGDS such as self-similar sets \cite{Trken}, \cite{Falconer_Minkowski}, \cite{Gatzouras}, \cite{Dundee}, \cite{Lapidus_Frankenhuysen_Springer}, \cite{LapPeaWin}, \cite{LapPom}, \cite{WinterLlorente}, \cite{Winter_thesis} and self-conformal sets \cite{Bohl}, \cite{Uta}, \cite{KesKom}, \cite{Diss}. New to this article is a general approach for the setting of cGDS which extends and recovers the previous results from the literature.
In this way, we can e.\,g.\ show that the Minkowski content of a limit set of a Fuchsian group of Schottky type always exists (see Sec.~\ref{sec:Fuchsian}), which proves a conjecture by M.\,L.\,Lapidus from 1993 that is stated in \cite{Dundee}.

A main motivation for studying the (local) Minkowski content and the (local) fractal Euler characteristic arises from fractal geometry, where one aims to find characteristics that describe the geometric structure of a fractal set. The (local) Minkowski content and the (local) fractal Euler characteristic can be viewed as such tools. They complement the notion of dimension and are capable of distinguishing between sets of the same Hausdorff- or Minkowski dimension. 
The Minkowski content of a set $Y\subset\mathbb R$ is the limit as $\eps$ tends to zero of the re-scaled length of the $\eps$-parallel neighbourhood of $Y$. Furthermore, for intervals, it coincides with the length of the interval. Therefore, the Minkowski content can be interpreted as ``fractal length''.
An interpretation of the fractal Euler characteristic is given by its name (see \cite{WinterLlorente}).
The local Minkowski content and the local fractal Euler characteristic are defined as weak limits of measures. They are Borel-measures which describe the ``fractal length'' and ``fractal Euler characteristic'' of a given fractal inside a Borel set.
When the weak limits exist, then their total masses respectively coincide with the Minkowski content and the fractal Euler characteristic.

In the literature, primarily the Minkowski content has received a lot of attention.
Especially Minkowski measurability (i.\,e.\ the existence of the Minkowski content in $(0,\infty)$) of self-similar subsets of $\mathbb R$ has been intensely studied. One important result is given for non-degenerate self-similar subsets of $\mathbb R$ whose associated iterated function system (IFS) consisting of similarities (sIFS) satisfies the open set condition with connected feasible open set. Such a set is Minkowski measurable if and only if the sIFS is non-lattice \cite{Falconer_Minkowski}, \cite{Dundee}, \cite{Lapidus_Frankenhuysen_Springer}, \cite{LapPom}.
We significantly extend this important result and provide an alternative proof by showing that the analogous statement is true also in the graph directed setting. To be more precise, we obtain that a non-degenerate limit set of a cGDS that consists of similarities (sGDS) is Minkowski measurable if and only if the sGDS is non-lattice (Cor.~\ref{cor:sGDSMink}). 
This convenient equivalence statement for systems consisting of similarities unfortunately fails to hold for general conformal systems: In \cite{KesKom} it has been shown that there exist non-degenerate self-conformal sets arising from lattice conformal IFS (cIFS) for which the Minkowski content and the fractal Euler characteristic exist. Since cIFS are special types of cGDS we cannot expect the equivalence to be valid for general limit sets of cGDS either. Indeed, in Prop.~\ref{thm:conformalMinkowski} we provide a sufficient condition under which the Minkowski content and the fractal Euler characteristic of a limit set of a lattice cGDS exist. However, one direction of the equivalence remains true: In the non-lattice case, the Minkowski content and the fractal Euler characteristic exist (see Rem.~\ref{rmk:Minkexistenceconf}). Moreover, in Rem.~\ref{rmk:Minkexistenceconf} we see that average versions of the Minkowski content and the fractal Euler characteristic always exist and provide explicit formulae to determine their values.

Let us now turn to the local quantities. These have been investigated in the context of fractal curvature measures in \cite{Winter_thesis}. 
If the ambient space is of dimension one, then there are two fractal curvature measures: The 0-th fractal curvature measure (which is the local fractal Euler characteristic) and the 1-st fractal curvature measure (which is the local Minkowski content). The term ``curvature'' is appropriate for higher dimensional ambient spaces but strictly speaking not in $\mathbb R$. Therefore, we will exclusively use the terms local Minkowski content and local fractal Euler characteristic in the present article.  
 We obtain that the local Minkowski content and the local fractal Euler characteristic exist for limit sets of non-lattice cGDS and are constant multiples of the associated $\mdim$-conformal measure, where $\mdim$ denotes the Minkowski dimension of the limit set (see Thm.~\ref{thm:curvatureresult}).
For limit sets of lattice sGDS we prove that these measures do not exist (see Thm.~\ref{thm:similars}). They neither exist for piecewise $\mathcal C^{1+\alpha}$-diffeomorphic images of limit sets of lattice sGDS (see Thm.~\ref{thm:imafcm}). This latter statement is important, since there are $\mathcal C^{1+\alpha}$-diffeomorphic images of limit sets of lattice sGDS for which the Minkowski content and the fractal Euler characteristic do exist (see Ex.~\ref{ex:thexample}).
Also for limit sets of lattice cGDS consisting of analytic maps, the local Minkowski content and the local fractal Euler characteristic do not exist (see Thm.~\ref{thm:curvatureresult}). However, we show that in the lattice situation average versions of the local Minkowski content and the local fractal Euler characteristic of a limit set of a cGDS always exist and are again constant multiples of the associated $\mdim$-conformal measure (see Thm.~\ref{thm:curvatureresult}). 

From the above results we deduce the following. The limit set of a Fuchsian group of Schottky type can be represented by a limit set of a non-lattice cGDS. Thus, as a consequence of Rem.~\ref{rmk:Minkexistenceconf}, its limit set is Minkowski measurable. As mentioned above, this result proves a conjecture by M.\,L.\,Lapidus from 1993 stated in \cite{Dundee}, which plays an important role in the context of the Weyl-Berry conjecture. The Weyl-Berry conjecture for fractal drums is a conjecture on the distribution of the eigenvalues of the Laplacian on a domain with a fractal boundary (see \cite{Falconer_Minkowski}, \cite{Dundee}, \cite{Lapidus_Frankenhuysen_Springer}, \cite{LapPom}).
It addresses the problem of describing `the relationship between the shape (geometry) of the drum and its sound (its spectrum).' \cite[p.\,1]{Lapidus_Frankenhuysen_Springer}
A more detailed exposition on the results from the literature and on the above mentioned conjecture will be given in Rem.~\ref{rmk:conjLapidus}.
Besides the motivation from fractal geometry, the Weyl-Berry conjecture is a main motivation for studying the Minkowski content, see e.\,g.\ \cite{survey}, \cite{Dundee}, \cite{Lapidus_Frankenhuysen_Springer}.
A third motivation arises from non-commutative geometry: In Connes' seminal book \cite{Connes} the notion of a non-commutative fractal geometry is developed. There, it is shown that the natural analogue of the volume of a compact smooth Riemannian spin$^c$ manifold for a fractal set in $\mathbb R$ is that of the Minkowski content. This idea is also reflected in \cite{Falconer_Samuel}, \cite{Guido_Isola}, \cite{Samuel}.

For an overview of the relevant literature and more background on the (local) Minkowski content and the (local) fractal Euler characteristic as well as an overview of the recent development of this research area, we refer the reader to \cite{KesKom} and the survey \cite{survey}.
Moreover, there are several recent articles concerning higher dimensional ambient spaces. In \cite{Gatzouras} it is shown that the Minkowski content of self-similar sets arising from non-lattice sIFS that satisfy the OSC exists. Alternative proofs of this result and further investigations on the lattice case are provided in \cite{Trken}, \cite{LapPeaWin}, where tube formulas and zeta-functions are used. Such tube formulas have been extended to limit sets of sGDS in \cite{Trken_graph-directed}. Minkowski measurability of self-conformal sets in higher dimensional ambient spaces has been studied in \cite{Bohl}, \cite{Diss}. There, it is shown that, under certain geometric conditions, a self-conformal set arising from a non-lattice cIFS is Minkowski measurable. Moreover, fractal curvature measures and their average versions are studied.

This article is organised as follows. In Sec.~\ref{sec:cGDS} we give the construction of cGDS and their limit sets. In Sec.~\ref{sec:results} we present our main results on the existence of the Minkowski content, the fractal Euler characteristic and their local versions. Sec.~\ref{sec:examples} is devoted to demonstrating how the new results can be applied to various classes of examples of limit sets of cGDS. Sec.~\ref{sec:preliminaries} to \ref{sec:specialcases} deal with the proofs of the main theorems. More precisely, in Sec.~\ref{sec:preliminaries} we provide some background and prove auxiliary results. With this preparation we provide the proofs of our main results concerning limit sets of cGDS (Thm.~\ref{thm:curvatureresult} and Prop.~\ref{thm:conformalMinkowski}) in Sec.~\ref{sec:proofs}. Sec.~\ref{sec:specialcases} is devoted to the proofs of Thms.~\ref{thm:similars} to \ref{thm:analytic} dealing with the special cases of sGDS as well as piecewise $\mathcal C^{1+\alpha}$-diffeomorphic images of limit sets of sGDS.

\section{Conformal Graph Directed Systems}\label{sec:cGDS}

A core text concerning conformal graph directed systems (cGDS) is \cite{Urbanski_Buch}. The class of cGDS generalises the notion of conformal iterated function systems and gives rise to a much richer class of fractal sets such as limit sets of Fuchsian groups. In Sec.~\ref{sec:examples} we give examples of classes of fractal sets which can be obtained via a cGDS. In this section, we present the relevant definitions.

\begin{definition}[Directed multigraph]
	A \emph{directed multigraph} $(V,E,i,t)$ consists of a finite set of vertices $V$, a finite set of directed edges $E$ and functions $i,t\colon E\to V$ which determine the initial and terminal vertex of an edge. The edge $e\in E$ goes from $i(e)$ to $t(e)$. Thus, the \emph{initial} and \emph{terminal vertices} of $e$ are $i(e)$ and $t(e)$ respectively.
\end{definition}

\begin{definition}[Incidence matrix]\label{defn:incidence}
	Given a directed multigraph $(V,E,i,t)$, an $(\card E)\times(\card E)$-matrix $A=(A_{e,e'})_{e,e'\in E}$ with entries in $\{0,1\}$, which satisfies $A_{e,e'}=1$ if and only if $t(e)=i(e')$ for edges $e,e'\in E$, is called an \emph{incidence matrix}. The incidence matrix $A$ is called \emph{aperiodic and irreducible} if there exists an $n\in \mathbb N$ such that the entries of the $n$-folded product $A^n$ are all positive.
\end{definition}

\begin{definition}[GDS]\label{def:GDS}
	A \emph{graph directed system (GDS)} consists of a directed multigraph $(V,E,i,t)$ with incidence matrix $A$, a family of non-empty compact connected metric spaces $(X_v)_{v\in V}$ and for each edge $e\in E$ an injective contraction $\phi_e\colon X_{t(e)}\to X_{i(e)}$ with Lipschitz constant less than or equal to $r$ for some $r\in (0,1)$. Briefly, the family $\Phi\defeq(\phi_e\colon X_{t(e)}\to X_{i(e)})_{e\in E}$ is called a GDS.
\end{definition}

In this paper, we consider fractal subsets of the real line. Therefore, we restrict the definition of a cGDS to the one-dimensional Euclidean space $(\mathbb R,\lvert\cdot\rvert)$.
For a subset $Y$ of $(\mathbb R,\lvert\cdot\rvert)$ we let $\inte(Y)$ denote its interior, $\overline{Y}$ its closure and $\partial Y\defeq \overline{Y}\setminus\inte(Y)$ its boundary.

\begin{definition}[cGDS]\label{def:cGDS}
	A GDS is called \emph{conformal (cGDS)} if 
	\begin{enumerate}
		\item for every $v\in V$, $X_v\subset\mathbb R$ is a compact interval with non-empty interior,
		\item the \emph{open set condition (OSC)} is satisfied, in the sense that, for distinct $e,e'\in E$ we have
		\[
			\phi_{e}(\inte(X_{t(e)}))\cap \phi_{e'}(\inte(X_{t(e')}))=\varnothing\qquad\text{and}
		\]
		\item for every vertex $v\in V$ there exists an open interval $W_v\supset X_v$ such that for every $e\in E$ with $t(e)=v$ the map $\phi_e$ extends to a $\mathcal C^{1+\alpha}$-diffeomorphism from $W_v$ into $W_{i(e)}$, whose derivative $\phi'_e$ is bounded away from zero on $W_v$, where $\alpha\in(0,1]$.
	\end{enumerate}
\end{definition}

We also consider the special case of cGDS where the contractions $\phi_e$ for $e\in E$ are similarities:

\begin{definition}[sGDS]
	A cGDS, whose maps $\phi_e$ are similarities for $e\in E$, is referred to as sGDS.
\end{definition}

\begin{remark}\label{rmk:GDMS}
  In the sequel, we will often refer to results from \cite{Urbanski_Buch}, where conformal graph directed Markov systems (cGDMS) are treated. Such systems differ from cGDS in that an incidence matrix for a cGDMS only fulfills the property that $A_{e,e'}=1$ implies $t(e)=i(e')$. 
  However, every cGDS is a cGDMS and conversely, a cGDMS in $\mathbb R$ can always be represented by a cGDS, namely by substituting $(\phi_e(X_{t(e)}))_{e\in E}$ in for the sets $(X_v)_{v\in V}$ and defining the edges accordingly. 
\end{remark}

In order to define the limit set of a cGDS, we fix a cGDS with the notation from Defs.~\ref{def:GDS} and \ref{def:cGDS}.
The set of \emph{infinite admissible words} given by the incidence matrix $A$ is defined to be
\begin{equation}\label{eq:EAinfty}
	E_A^{\infty}\defeq \{\om=\om_1\om_2\cdots\in E^{\mathbb N}\mid A_{\om_n,\om_{n+1}}=1\ \text{for all}\ n\in \mathbb N\}.
\end{equation}
The set of sub-words of length $n\in\mathbb N$ is denoted by $E_A^n$ and the set of all finite sub-words including the empty word $\varnothing$ by $E_A^*$.
For a finite word $\om\in E_A^*$ we let $n(\om)$ denote its length, where $n(\varnothing)\defeq 0$, define $\phi_{\varnothing}$ to be the identity map on $\bigcup_{v\in V}X_v$ and for $\om\in E_A^*\setminus\{\varnothing\}$ set
\[
	\phi_{\om}\defeq\phi_{\om_1}\circ\cdots\circ\phi_{\om_{n(\om)}}\colon X_{t(\om_{n(\om)})}\to X_{i(\om_1)},
\]
where we let $\om_i$ denote the $i$-th letter of the word $\om$ for $i\in\{1,\ldots,n(\om)\}$, i.\,e.\  $\om=\om_1\cdots\om_{n(\om)}$.
For two finite words $\omneu=\omneu_1\cdots\omneu_n$, $\om=\om_1\cdots\om_m\in E_A^{*}$ with $A_{\omneu_n,\om_1}=1$, we let $\omneu\om\defeq\omneu_1\cdots\omneu_n\om_1\cdots\om_m\in E_A^*$ denote their concatenation. Likewise, we set $\omneu\om\defeq\omneu_1\cdots\omneu_n\om_1\om_2\cdots$ for $\omneu=\omneu_1\cdots\omneu_n\in E_A^*$ and $\om=\om_1\om_2\cdots\in E_A^{\infty}$ with $A_{\omneu_n,\om_1}=1$.
For an infinite word $\om=\om_1\om_2\cdots \in E_A^{\infty}$ and $n\in \mathbb N$ the \emph{initial word of length $n$} is defined to be $\om\vert_{n}\defeq \om_1\cdots\om_n$. 

For $\om\in E_A^{\infty}$ the sequence $(\phi_{\om\vert_n}(X_{t(\om_n)}))_{n\in\mathbb N}$ is a descending sequence of non-empty compact sets and therefore $\bigcap_{n\in \mathbb N} \phi_{\om\vert_n}(X_{t(\om_n)})\neq\varnothing$. Recall from Def.~\ref{def:GDS} that $r\in(0,1)$ denotes a common Lipschitz constant of $\phi_e$ for $e\in E$. Since $\diam(\phi_{\om\vert_n}(X_{t(\om_n)}))\leq r^n\diam(X_{t(\om_n)})\leq r^n\max\{\diam(X_v)\mid v\in V\}$  for every $n\in\mathbb N$, the intersection 
\[
	\bigcap_{n\in\mathbb N}\phi_{\om\vert_n}(X_{t(\om_n)})
\]
is a singleton and we denote its only element by $\bij(\om)$. The \emph{code map} is defined to be the map $\bij\colon E_A^{\infty}\to\bigcup_{v\in V} X_v$ given by $\om\mapsto\bij(\om)$. 

\begin{definition}[Limit set of a cGDS]
	The \emph{limit set} of the cGDS $(\phi_e)_{e\in E}$ is defined to be 
	\[
		F\defeq\bij(E_A^{\infty}).
	\]
\end{definition}
Limit sets of cGDS often have a fractal structure. They include invariant sets of conformal iterated function systems, the so-called self-conformal sets, as well as self-similar sets. These are defined as follows.
\begin{definition}[cIFS, sIFS, self-conformal set, self-similar set]\label{defn:CIFS}
	A \emph{conformal iterated function system (cIFS)} is a cGDS $\Psi\defeq(\psi_1,\ldots,\psi_N)$ whose set of vertices $V$ is a singleton and whose set of edges contains at least two elements.
	The unique limit set of a cIFS is called the \emph{self-conformal set} associated with $\Psi$.
	In the case that the maps $\psi_1,\ldots,\psi_N$ are similarities, the limit set is called the \emph{self-similar set} associated with $\Psi$ and $\Psi$ is referred to as an \emph{sIFS}.
\end{definition}

 In order to show the significance of cGDS, Sec.~\ref{sec:examples} is devoted to examples of important classes of such sets.

\section{Main Results}\label{sec:results}

\subsection{Notation, Definitions and First Results}\label{sec:Notation}
Before stating our results, let us begin with recalling the relevant notation and definitions, in particular the local Minkowski content and the local fractal Euler characteristic. For further background we refer the reader to \cite{KesKom}.

We let $\leb^0$ and $\leb^1$ respectively denote the counting measure and the one-dimensional Lebesgue measure. For an arbitrary subset $Y\subset\mathbb R$ and $\eps>0$ we define $Y_{\eps}\defeq\{x\in\mathbb R\mid\inf_{y\in Y}\lvert x-y\rvert\leq\eps\}$ to be the \emph{$\eps$-parallel neighbourhood} of $Y$. For the remainder of this section we assume that $Y$ is non-empty and compact. The \emph{1-st} and \emph{0-th scaling exponents} of $Y$ are respectively defined to be
$s_1(Y)\defeq \inf\{t\in\mathbb R\mid \eps^t\leb^1(Y_{\eps})\to 0\ \text{as}\ \eps\to 0\}$ and 
$s_0(Y)\defeq \inf\{t\in\mathbb R\mid \eps^t\leb^0(\partial Y_{\eps})\to 0\ \text{as}\ \eps\to 0\}$.

\begin{definition}\label{defn:curvaturemeasures}\label{def:averagefc}
	Provided, that the weak limit
	\[
	C_1^f(Y,\cdot)\defeq\wlim_{\eps\to 0}\eps^{s_1(Y)}\leb^1(Y_{\eps}\cap\cdot)
	\]
	of the finite Borel measures $\eps^{s_1(Y)}\leb^1(Y_{\eps}\cap\cdot)$ exists, we call $C_1^f(Y,\cdot)$ the \emph{local Minkowski content} of $Y$. Likewise, the weak limit
	\[
	C^f_0(Y,\cdot)\defeq\wlim_{\eps\to 0} \eps^{s_0(Y)}\leb^0(\partial Y_{\eps}\cap\cdot)/2
	\]
	is called the \emph{local fractal Euler characteristic} of $Y$, if it exists. 
	Moreover, provided that the weak limits exist, we respectively call
	\begin{align*}
	\widetilde{C}_1^f(Y,\cdot)&\defeq\wlim_{T\searrow 0} \lvert\ln T\rvert^{-1}\int_T^1\eps^{s_1(Y)-1}\leb^1(Y_{\eps}\cap\cdot)\textup{d}\eps\qquad\text{and}\\
	\widetilde{C}_0^f(Y,\cdot)&\defeq\wlim_{T\searrow 0}\lvert\ln T\rvert^{-1}\int_T^1\eps^{s_0(Y)-1}\leb^0(\partial Y_{\eps}\cap\cdot)\textup{d}\eps/2
	\end{align*}
	the \emph{average local Minkowski content} of $Y$ and the \emph{average local fractal Euler characteristic} of $Y$.
	Moreover, for a Borel set $B\subseteq\mathbb R$ we set 
  \begin{equation*}\begin{array}[h]{rlrl}
    \overline{C}_0^f(Y,B)
    &\hspace{-.3cm}\defeq\displaystyle{\limsup_{\eps\to 0}\eps^{s_0(Y)}\leb^0(\partial Y_{\eps}\cap B)/2},
    &\overline{C}_1^f(Y,B)&\hspace{-.3cm}\defeq\displaystyle{\limsup_{\eps\to 0}\eps^{s_1(Y)}\leb^1(Y_{\eps}\cap B)}\\
    \underline{C}_0^f(Y,B)
    &\hspace{-.3cm}\defeq\displaystyle{\liminf_{\eps\to 0}\eps^{s_0(Y)}\leb^0(\partial Y_{\eps}\cap B)/2},
    &\underline{C}_1^f(Y,B)&\hspace{-.3cm}\defeq\displaystyle{\liminf_{\eps\to 0}\eps^{s_1(Y)}\leb^1(Y_{\eps}\cap B)}.
  \end{array}\end{equation*}
\end{definition}

Notice, if $C_k^f(Y,\cdot)$ exists, then also $\widetilde{C}_k^f(Y,\cdot)$ exists and the two Borel measures coincide. However, we will see that $C_k^f(Y,\cdot)$ does not always exist, whereas we will prove that the average version $\widetilde{C}_k^f(Y,\cdot)$ always exists for limit sets of cGDS.

\begin{remark}
  The (average) fractal Euler characteristic was investigated for self-similar sets in \cite{WinterLlorente}. In higher dimensional ambient spaces, the local Minkowski content and the local fractal Euler characteristic belong to the class of fractal curvature measures as introduced by S.\,Winter in \cite{Winter_thesis}. However, the notion of curvature is appropriate only in higher dimensional ambient spaces and therefore not in the context of the present article.
\end{remark}

The first step towards determining the quantities from Def.~\ref{defn:curvaturemeasures} for limit sets of cGDS is to evaluate the scaling exponents $s_1$ and $s_0$. Scaling exponents have been intensely studied in \cite{PokornyWinter}. For a limit set $F$ of a cGDS with aperiodic and irreducible incidence matrix, $s_1(F)$ and $s_0(F)$ are strongly linked to the Minkowski dimension $\textup{dim}_M(F)$ of $F$ which is proven to exist for such sets $F$ (see Prop.~\ref{thBedford}). For an arbitrary bounded set $Z\subset\mathbb R$ the Minkowski dimension is defined by 
\begin{equation}\label{eq:Minkdimension}
  \textup{dim}_{M}(Z)\defeq 1-\lim_{\eps\searrow 0}\frac{\ln\leb^1(Z_{\eps})}{\ln\eps},
\end{equation}
whenever this limit exists. When $\textup{dim}_{M}(Z)$ exists then also the box-counting dimension of $Z$ exists and both quantities coincide (see \cite[Prop.~3.2]{Falconer_Foundation}). The connection between the Minkowski dimension and the scaling exponents is provided in the next proposition.

\begin{proposition}\label{prop:interval}\label{prop:intervalfractal}
  Let $\mdim$ denote the Minkowski dimension of the limit set $F$ of a cGDS $\Phi$ with aperiodic and irreducible incidence matrix. Then either $\leb^1(F)=0$, in which case $s_1(F)=\mdim-1$ and $s_0(F)=\mdim$, or $F$ is a finite union of compact intervals with non-empty interior. For any finite union $Z\subset\mathbb R$ of compact intervals with non-empty interior we have $s_1(Z)=s_0(Z)=0$.
\end{proposition}
\begin{proof}
  First, we show that either $\leb^1(F)=0$ or $F$ is a finite union of compact intervals with non-empty interior.
  For $n\in\mathbb N$ define $X^{(n)}\defeq\bigcup_{\om\in E_A^{n}}\phi_{\om}(X_{t(\om_n)})$ and set $X\defeq\bigcup_{v\in V}X_v$. If $\leb^1(\inte(X)\setminus X^{(1)})>0$, then $\leb^1(F)=0$ by \cite[Prop.~4.5.9]{Urbanski_Buch}. On the other hand, if $\leb^1(\inte(X)\setminus X^{(1)})=0$, then $\leb^1(X\setminus X^{(1)})=0$, since the cardinality of $\partial X$ is finite. It follows that $X\setminus X^{(1)}=\varnothing$, as both $X$ and $X^{(1)}$ are finite unions of compact intervals. Clearly then $X^{(n)}=X$ for all $n\in\mathbb N$ and $F=X$.
  Next, we turn to the connection between $\mdim$ and the scaling exponents.
  If $\lambda^1(F)=0$, the equality $s_1(F)=\mdim-1$ follows straight from the definitions of $s_1(F)$ and $\mdim$. Moreover, the relation $s_0(F)=\mdim$ is a consequence of $s_1(F)=\mdim-1$ and \cite[Cor.~3.2]{RatajWinter}.
  If $Z$ is a finite union of compact intervals with non-empty interior, then $\lambda^1(Z)>0$ and $\lambda^0(\partial Z_{\eps})$ is positive and uniformly bounded in $\eps$, which imply $s_1(Z)=s_0(Z)=0$.
\end{proof}

In the setting of Prop.~\ref{prop:interval}, if $\lambda^1(F)=0$ then we call $F$ \emph{non-degenerate}; otherwise we call $F$ \emph{degenerate}.

\begin{remark}In the degenerate situation an immediate consequence of $s_1(Z)=s_0(Z)=0$ from the above proposition is that both the local Minkowski content and the local fractal Euler characteristic of $Z$ exist and satisfy
\begin{equation*}
	C^f_1(Z,\cdot)=\leb^1(Z\cap\cdot)\quad \text{and}\quad
	C^f_0(Z,\cdot)=\leb^0(\partial Z\cap\cdot)/2 .
\end{equation*}
\end{remark}
The more interesting case of Prop.~\ref{prop:interval} is the non-degenerate case. 
For stating our results for such sets, we fix a cGDS $(V,E,i,t,A)$ and assume that the incidence matrix $A$ is aperiodic and irreducible (see Def.~\ref{defn:incidence}). Let $(X_v)_{v\in V}$ denote the associated compact intervals with non-empty interior and let $\Phi\defeq(\phi_e\colon X_{t(e)}\to X_{i(e)})_{e\in E}$ denote the family of injective $r$-Lipschitz maps for some $r\in(0,1)$. Further, let $F$ denote the unique limit set and let $\mdim\defeq\dim_M(F)$ be its Minkowski dimension. 
A central role with regard to our results is played by the geometric potential function:
\begin{definition}[Geometric potential function, shift-map]\label{defn:geompotfcn}
The \emph{geometric potential function} $\xi\colon E_A^{\infty}\to\mathbb R$ is defined by $\xi(\om)\defeq -\ln\lvert\phi'_{\om_1}(\bij(\sigma\om))\rvert$ for $\om=\om_1\om_2\cdots\in E_A^{\infty}$. 
Here $\sigma\colon E_A^*\cup E_A^{\infty}\to E_A^*\cup E_A^{\infty}$ denotes the \emph{shift-map} which is defined by 
$\sigma(\om_1\om_2\cdots)\defeq\om_2\om_3\cdots\in E_A^{\infty}$ for $\om_1\om_2\cdots\in E_A^{\infty}$,
$\sigma(\om_1\cdots\om_n)\defeq\om_2\cdots\om_n\in E_A^{n-1}$ for $\om_1\cdots\om_n\in E_A^n$, where $n\geq 2$ and
$\sigma(\om)\defeq \varnothing$ for $\om\in\{\varnothing\}\cup E_A^1$.
\end{definition}

We equip $E^{\mathbb N}$ with the product topology of the discrete topologies on $E$ and equip the set of infinite admissible words $E_A^{\infty}\subset E^{\mathbb N}$ with the subspace topology. This is the weakest topology with respect to which the canonical projections onto the coordinates are continuous. 
The space of continuous real-valued functions on $E_A^{\infty}$ is denoted by $\mathcal C (E_A^{\infty})$.
Note that the geometric potential function $\xi$ belongs to $\mathcal C(E_A^{\infty})$. A crucial property of $\xi$ is whether it is lattice or non-lattice.
\begin{definition}[Co-homologous, lattice, non-lattice]\label{def:nonlattice}\label{def:nonlatticefractal}~\\[-2.5ex]
  \begin{enumerate}
  \item Functions $f_1,$ $f_2\in\mathcal C( E_A^{\infty})$ are called \emph{co-homologous}, if there exists a function $\psi\in\mathcal C( E_A^{\infty})$ such that $f_1-f_2=\psi-\psi\circ\sigma$. A function $f\in\mathcal C( E_A^{\infty})$ is said to be \emph{lattice}, if $f$ is co-homologous to a function whose range is contained in a discrete subgroup of $\mathbb R$. Otherwise, we say that $f$ is \emph{non-lattice}.
  \item If the geometric potential function $\xi$ is non-lattice, then we call the cGDS $\Phi$ \emph{non-lattice}. On the other hand, if $\xi$ is lattice, then we call $\Phi$ \emph{lattice}.
  \end{enumerate}
\end{definition}

We let $\entro(\mu_{-\mdim\xi})$ denote the measure theoretical entropy of the shift-map $\sigma$ with respect to the unique $\sigma$-invariant Gibbs measure $\mu_{-\mdim\xi}$ for the potential function $-\mdim\xi$ (see \eqref{eq:entropy} for a definition).
The unique probability measure $\nu$ supported on $F$, which for all distinct $e,e'\in E$ satisfies
\begin{equation}\label{eq:conformalmeasure}
	\nu(\phi_e (X_{t(e)})\cap\phi_{e'}(X_{t(e')}))=0 \quad\text{and}\quad	\nu(\phi_{e}B)=\int_B \lvert\phi'_{e}\rvert^{\mdim}\textup{d}\nu
\end{equation}
for all Borel sets $B\subseteq X_{t(e)}$ is called the \emph{$\mdim$-conformal measure} associated with $\Phi$. Uniqueness and existence is provided in \cite[Thm.~4.2.9]{Urbanski_Buch} and goes back to the work of \cite{DenkerUrbanski}, \cite{Patterson}, \cite{Sullivan}.

For a vertex $v\in V$ we denote the set of edges whose initial and respectively terminal vertex is $v$ by
\[
	I_v\defeq\{e\in E\mid i(e)=v\}\quad\text{and}\quad 
	T_v\defeq\{e\in E\mid t(e)=v\}.
\]
Moreover, for $n\in\mathbb N$ we set
\begin{equation*}
\begin{array}{rlrl}
	I_v^n&\defeq\ \{\om\in E_A^n\mid i(\om_1)=v\},\qquad
	&T_v^n&\defeq\ \{\om\in E_A^n\mid t(\om_n)=v\},\\
	I_v^*&\defeq\ \bigcup_{k\in\mathbb N} I_v^k,
	&T_v^*&\defeq\ \bigcup_{k\in\mathbb N} T_v^k\quad \text{and}\quad\\
	I_v^{\infty}&\defeq\ \{\om\in E_A^{\infty}\mid i(\om_1)=v\}.&&
\end{array}
\end{equation*}
For a finite word $\om\in E_A^*$ the \emph{$\om$-cylinder set} is defined to be 
\[
	[\om]\defeq\{\omneu\in E_A^{\infty}\mid \omneu_i=\om_i\ \text{for}\ i\in\{1,\ldots,n(\om)\}\},\ \text{in particular}\ [\varnothing]=E_A^{\infty}.
\]

Fundamentally important objects in our main statements are the primary gaps of $F$ and their images. These are certain intervals in the complement of the limit set, which are defined in the following way. Set
\begin{equation}\label{eq:Lv}
	L^v\defeq\left\langle\bigcup_{e\in I_v}\bij[e]\right\rangle\setminus \bigcup_{e\in I_v}\left\langle\bij[e]\right\rangle,
\end{equation}
where $v\in V$ and $\langle Y\rangle$ denotes the convex hull of $Y$.
We let $n_v$ denote the number of connected components of $L^v$. In Prop.~\ref{prop:Lvnonempty} we show that $\bigcup_{v\in V} L^v\neq\varnothing$ if $\leb^1(F)=0$, hence, $\sum_{v\in V}n_v\geq 1$.
If $L^v\neq\varnothing$, we denote the connected components of $L^v$ by $L^{v,j}$, where $j$ ranges over $\{1,\ldots,n_v\}$ and call the sets $L^{v,j}$ the \emph{primary gaps} of $F$. For every $\om\in T_v^*$ we define $L_{\om}^{v,j}\defeq\phi_{\om}(L^{v,j})$ and call these sets the \emph{\main\ gaps} of $F$.

\subsection{Exposition of the Main Results}\label{sec:mainresults}
Now, we are able to present our main results and for this purpose fix the notation from Sec.~\ref{sec:Notation}. In particular, let $\Phi\defeq(\phi_e)_{e\in E}$ denote a cGDS with aperiodic and irreducible incidence matrix and let $F$ denote its limit set. Set $\mdim\defeq\dim_M(F)$ and let $\xi$ denote the geometric potential function associated with $\Phi$.
Further, denote by $\entro(\mu_{-\mdim\xi})$ the measure theoretical entropy of the shift-map $\sigma$ with respect to the unique shift-invariant Gibbs measure $\mu_{-\mdim\xi}$ for the potential function $-\mdim\xi$ (see Sec.~\ref{sec:PF}). 
The proofs of the theorems of this subsection are presented in Sec.~\ref{sec:proofs} and \ref{sec:specialcases}.

\begin{theorem}\label{thm:curvatureresult}
  Assume that $\leb^1(F)=0$. Then the following hold.
  \begin{enumerate}
  \item\label{curvatureresult:average} The average local Minkowski content and the average local fractal Euler characteristic of $F$ always exist and are constant multiples of the $\mdim$-conformal measure $\nu$ associated with $\Phi$, i.\,e.
	\[
	\hspace{0.5cm}\widetilde{C}_1^f(F,\cdot)=\frac{2^{1-\mdim}c}{(1-\mdim)\entro(\mu_{-\mdim\xi})}\cdot\nu(\cdot)\quad\text{and}\quad
	\widetilde{C}_0^f(F,\cdot)=\frac{2^{-\mdim}c}{\entro(\mu_{-\mdim\xi})}\cdot\nu(\cdot), 
	\]
    where the constant $c$ is given by the well-defined positive and finite limit
    \begin{equation}\label{eq:constantthm}
      c\defeq \lim_{m\to\infty}\sum_{v\in V}\sum_{j=1}^{n_v}\sum_{\om\in T_v^m} \lvert L_{\om}^{v,j}\rvert^{\mdim}.
    \end{equation}
  \item\label{curvatureresult:nonlattice} If $\xi$ is non-lattice, then both the local Minkowski content and the local fractal Euler characteristic of $F$ exist and satisfy $C_k^f(F,\cdot)=\widetilde{C}_k^f(F,\cdot)$ for $k\in\{0,1\}$.
  \item\label{curvatureresult:lattice} If $\xi$ is lattice, then there exists a constant $\overline{c}\in\mathbb R$ such that 
    $0<\underline{C}^f_k(F,\mathbb R)\leq\overline{C}^f_k(F,\mathbb R)\leq \overline{c}$ for $k\in\{0,1\}$. If additionally the system $\Phi$ consists of analytic maps, then neither the local Minkowski content nor the local fractal Euler characteristic of $F$ exists. 
  \end{enumerate}
\end{theorem}
Thm.~\ref{thm:curvatureresult}\ref{curvatureresult:nonlattice} and \ref{curvatureresult:lattice} in particular show that the scaling exponents of $F$ can alternatively be characterised by $s_1(F)=\sup\{t\in\mathbb R\mid \eps^t\leb^1(F_{\eps})\to\infty\ \text{as}\ \eps\to 0\}$ and $s_0(F)= \sup\{t\in\mathbb R\mid \eps^t\leb^0(\partial F_{\eps})\to\infty\ \text{as}\ \eps\to 0\}$. 
Since $\overline{C}^f_k(F,\cdot)$ is monotonically increasing as a set function in the second component, Thm.~\ref{thm:curvatureresult}\ref{curvatureresult:lattice} also shows that $\overline{C}^f_k(F,B)\leq\overline{c}$ for all Borel sets $B\subset\mathbb R$.

\begin{theorem}[sGDS]\label{thm:similars}
  Suppose that $\Phi$ is an sGDS.
  Assume that $\leb^1(F)=0$ and let $\eigenf_{-\mdim\xi}$ denote the unique strictly positive eigenfunction with eigenvalue one of the Perron-Frobenius operator for the potential function $-\mdim\xi$ (see Sec.~\ref{sec:PF}).
  Then, additionally to the statements of Thm.~\ref{thm:curvatureresult}, the following hold. 
  \begin{enumerate}
  \item\label{ss:average} 
    The constant $c$ from \eqref{eq:constantthm} simplifies to the finite sum
    \begin{align*}
      c=\sum_{v\in V}\sum_{j=1}^{n_v}\eigenf_{-\mdim\xi}(\om^v)\lvert L^{v,j}\rvert^{\mdim},
    \end{align*}
    which is independent of the choice of $\om^v\in I_v^{\infty}$.
  \item\label{ss:lattice} If $\xi$ is lattice, then the following holds. For $k\in\{0,1\}$ and for every Borel set $B\subseteq\mathbb R$ for which $F\cap B$ is a non-empty finite union of sets of the form $\bij[\om]$, where $\om\in E_A^*$, and for which $F_{\eps}\cap B=(F\cap B)_{\eps}$ for all sufficiently small $\eps>0$ we have that
    \begin{equation*}
      0<\underline{C}_k^f(F,B)<\overline{C}_k^f(F,B)<\infty.
    \end{equation*}
    Consequently neither the local Min\-kowski content nor the local fractal Euler characteristic of $F$ exists.
  \end{enumerate}
\end{theorem}

An interesting subclass of limit sets of cGDS is the class of piecewise $\mathcal C^{1+\alpha}$-diffeomorphic images of limit sets of sGDS, where $\alpha\in(0,1]$ and $\mathcal C^{1+\alpha}$ denotes the class of real-valued functions which are differentiable with $\alpha$-H\"older continuous derivative. A nice relationship between the (average) local Minkowski content and the (average) local fractal Euler characteristic of the limit set of the sGDS and of its piecewise $\mathcal C^{1+\alpha}$-diffeomorphic image is provided in the next theorem. 
The analogue statements of Thm.~\ref{thm:imafcm}\ref{imafcm_average} and \ref{imafcm_nonlattice} have been obtained in \cite{Uta} for conformal $\mathcal{C}^{1+\alpha}$-diffeomorphic images of self-similar sets in higher dimensional ambient spaces.

\begin{theorem}[Piecewise $\mathcal C^{1+\alpha}$-diffeomorphic images of limit sets of sGDS]\label{thm:imafcm}
  Let $R$ denote an sGDS with aperiodic and irreducible incidence matrix, with associated directed multigraph $(V,E,i,t)$ and with associated compact non-empty intervals $(Y_v)_{v\in V}$. Let $\sset\subset\mathbb R$ denote the limit set of $R$ and assume that $\leb^1(\sset)=0$. For each $v\in V$ let $g_v\colon W_v\to\mathbb R$ denote a $\mathcal C^{1+\alpha}(W_v)$-diffeomorphism which is defined on a connected open neighbourhood $W_v\subset\mathbb R$ of $Y_v$ such that $\lvert g_v'\rvert$ is bounded away from zero on $W_v$ and such that the interiors of $X_v\defeq g_v(Y_v)$ are pairwise disjoint and $\alpha\in(0,1]$. Set $F\defeq\bigcup_{v\in V} g_v(\sset\cap Y_v)$. Then we have the following.
    \begin{enumerate}
    \item\label{imafcm_average} The average local Minkowski content and the average local fractal Euler characteristic of both $\sset$ and $F$ exist. Moreover, $\widetilde{C}_k^f(F,\cdot)$ is absolutely continuous with respect to the push-forward measure $\widetilde{C}_k^f(\sset,\bigcup_{v\in V}g_v^{-1}(\cdot))$ for $k\in\{0,1\}$. Their Radon-Nikodym derivative is, for $v\in V$ and $k\in\{0,1\}$, given by
      \begin{equation*}
	\frac{\textup{d}\widetilde{C}_k^f(F,\cdot)}{\textup{d}\widetilde{C}_k^f(\sset,\bigcup_{v'\in V}g_{v'}^{-1}(\cdot))}\Bigg{\rvert}_{X_v}=\lvert g_v'\circ g_v^{-1}\rvert^{\mdim}\Bigg{\rvert}_{X_v},
      \end{equation*}
      where $\mdim\defeq\dim_M(\sset)$ denotes the Minkowski dimension of $\sset$.
    \item\label{imafcm_nonlattice} If $R$ is non-lattice, then the local Minkowski content and the local fractal Euler characteristic of both $\sset$ and $F$ exist and coincide with the respective average versions.
    \item\label{imafcm_lattice} If $R$ is lattice, then neither the local Minkowski content nor the local fractal Euler characteristic of $\sset$ and $F$ exist.
    \end{enumerate}
\end{theorem}

Piecewise $\mathcal C^{1+\alpha}$-diffeomorphic images of limit sets of sGDS play an important role in the theory of general lattice cGDS. Namely, if a lattice cGDS consists of analytic maps, then its limit set $F$ is an image of a limit set of an sGDS under a piecewise $\mathcal C^{1+\alpha}$-diffeomorphism:

\begin{theorem}[Rigidity]\label{thm:analytic}
	Let $\Phi$ be a lattice cGDS consisting of analytic maps and let $F\subset\mathbb R$ denote its limit set. Then there exists a limit set $\sset\subset\mathbb R$ of a lattice sGDS, with associated compact intervals $(Y_v)_{v\in V}$ and analytic maps $g_v\colon W_v\to\mathbb R$ with $\lvert g_v'\rvert$ bounded away from zero, where $W_v$ is an open neighbourhood of $Y_v$, such that $F=\bigcup_{v\in V}g_v (K\cap Y_v)$.
\end{theorem}

Thm.~\ref{thm:analytic} is a generalisation of \cite[Thm.~2.2]{KesKom}, which addresses cIFS.

\subsection{Results on the Minkowski Content}
The theorems from the preceding subsection immediately imply results on the existence and value of the (average) Minkowski content of limit sets of cGDS. 
The (average) Minkowski content is defined as follows. Let $Y\subset\mathbb R$ denote a set whose Minkowski dimension $\textup{dim}_M(Y)\eqdef\mdim$ exists.
The \emph{upper Minkowski content} $\overline{\mathcal M}(Y)$ and the \emph{lower Minkowski content} $\underline{\mathcal M}(Y)$ of $Y$ are respectively defined to be
\begin{equation}
  \overline{\mathcal M}(Y)
  \defeq\limsup_{\eps\to 0}\eps^{\mdim-1}\leb^1(Y_{\eps})\quad \text{and}\quad
  \underline{\mathcal M}(Y)
  \defeq\liminf_{\eps\to 0}\eps^{\mdim-1}\leb^1(Y_{\eps}).
\end{equation}
If the upper and lower Minkowski contents coincide, then we denote the common value by $\mathcal M(Y)$ and call it the \emph{Minkowski content} of $Y$. In the case that the Minkowski content exists, is positive and finite, we call $Y$ \emph{Minkowski measurable}. The \emph{average Minkowski content} of $Y$ is defined to be the following limit, provided it exists
\begin{equation}
  \widetilde{\mathcal M}(Y)\defeq\lim_{T\searrow 0}\lvert\ln T\rvert^{-1}\int_T^1\eps^{\mdim-2}\leb^1(Y_{\eps})\textup{d}\eps.
\end{equation}
Analogously, one defines the \emph{fractal Euler characteristic} of $Y$ by 
\begin{equation}
	C_0^f(Y)\defeq\lim_{\eps\to 0}\eps^{s_0(Y)}\lambda^0(\partial F_{\eps})/2,
\end{equation}
whenever this limit exists.

We use the notation from the beginning of Sec.~\ref{sec:mainresults} for stating the implications of Thms.~\ref{thm:curvatureresult} to \ref{thm:imafcm} regarding the existence and the value of the Minkowski content.
\begin{remark}\label{rmk:Minkexistenceconf}
  Suppose that $\lambda^1(F)=0$. Immediate consequences of Prop.~\ref{prop:interval} and Thm.~\ref{thm:curvatureresult} are that 
  \begin{equation}\label{eq:Mcaverage}
    \widetilde{\mathcal{M}}(F)
    =\widetilde{C}_1^f(F,\mathbb R)
    =\frac{2^{1-\mdim}c}{(1-\mdim)\entro(\mu_{-\mdim\xi})}
  \end{equation}
  with $c$ as in \eqref{eq:constantthm}, and if $\xi$ is non-latice that $\mathcal{M}(F)$ exists and satisfies
  \begin{equation}\label{eq:Mcnonlattice}
    \mathcal{M}(F)
    =C_1^f(F,\mathbb R)
    =\widetilde{C}_1^f(F,\mathbb R)
    =\widetilde{\mathcal{M}}(F).
  \end{equation}
\end{remark}
In the above remark we addressed the average Minkowski content and the non-lattice case. The remaining lattice case is delicate with regard to Minkowski measurability, namely both existence and non-existence of the Minkowski content is possible.
A sufficient condition under which the Minkowski content exists is given in the following proposition, which will be proved in Sec.~\ref{sec:proofs}. Here, for an $\alpha$-H\"older continuous function $f\in\mathcal F_{\alpha}(E_A^{\infty})$ (see Sec.~\ref{sec:PF}) we let $\nu_f$ denote the unique eigenmeasure with eigenvalue 1 of the dual of the Perron-Fro\-benius operator for the potential function $f$ (see Sec.~\ref{sec:PF}).

\begin{proposition}\label{thm:conformalMinkowski}
  Assume that $\leb^1(F)=0$ and that $\xi$ is lattice. Then we have
  \begin{equation}\label{eq:Mclattice}
    0<\underline{\mathcal M}(F)\leq\overline{\mathcal M}(F)<\infty.
  \end{equation}
  Further, equality in \eqref{eq:Mclattice} can be attained. More precisely, let $\ze,\psi\in\mathcal{C}( E_A^{\infty})$ denote two functions satisfying $\xi-\ze=\psi-\psi\circ\sigma$, where the range of $\ze$ is contained in a discrete subgroup of $\mathbb R$ and $\aaa\in\mathbb R$ is maximal such that $\ze( E_A^{\infty})\subseteq\aaa\mathbb Z$. If, for every $t\in[0,\aaa)$, we have that
    \begin{align}\label{eq:existencecondition}
      &\sum_{n\in\mathbb Z}\ee^{-\mdim\aaa n}\nu_{-\mdim\ze}\circ\psi^{-1}([n\aaa,n\aaa+t))\nonumber\\
	&\qquad=\frac{\ee^{\mdim t}-1}{\ee^{\mdim\aaa}-1}\sum_{n\in\mathbb Z}\ee^{-\mdim\aaa n}\nu_{-\mdim\ze}\circ\psi^{-1}([n\aaa,(n+1)\aaa)),
    \end{align}
    then $\underline{\mathcal{M}}(F)=\overline{\mathcal{M}}(F)$, where the sums occurring in \eqref{eq:existencecondition} are finite sums.
\end{proposition}

Condition \eqref{eq:existencecondition} was obtained in \cite[(2.3)]{KesKom} as a condition implying that $\underline{\mathcal M}(F)=\overline{\mathcal M}(F)$ for self-conformal sets $F$ arising from a lattice cIFS. Thus, the above proposition states that the exact same condition implies Minkowski measurability for limit sets of cGDS. The necessary adaptations of the proof of \cite{KesKom} are outlined in Sec.~\ref{sec:proofs}.
An example of a lattice limit set of a cGDS, which satisfies \eqref{eq:existencecondition} and thus is Minkowski measurable, is given in Ex.~\ref{ex:thexample}.
However, in the special case, when the maps $\phi_e$ of the lattice cGDS are similarities, \eqref{eq:existencecondition} cannot be satisfied which  follows from Thm.~\ref{thm:similars}\ref{ss:lattice}. Indeed, a consequence of Thms.~\ref{thm:curvatureresult} and \ref{thm:similars} is the following.

\begin{corollary}\label{cor:sGDSMink}
  Suppose that $\Phi$ is an sGDS and that $\lambda^1(F)=0$. Then $F$ is Minkowski measurable if and only if the sGDS $\Phi$ is non-lattice.
\end{corollary}
This corollary provides an important extension of the result for sIFS given in \cite{Falconer_Minkowski}, \cite{Dundee}, \cite{Lapidus_Frankenhuysen_Springer},  \cite{LapPom} (see Rem.~\ref{rmk:conjLapidus}\ref{it:conjss}).

\begin{remark}
  Combining \eqref{eq:Mcaverage} with Thm.~\ref{thm:similars}\ref{ss:average}  one obtains an explicit expression for $\widetilde{\mathcal M}(F)$, when $\Phi$ is an sGDS.
\end{remark}

We can now turn to the class of piecewise $\mathcal C^{1+\alpha}$-diffeo\-mor\-phic images of limit sets arising from sGDS. Also here, the lattice case is delicate and so we first treat the average and non-lattice situations.

\begin{remark}
  Suppose that we are in the situation of Thm.~\ref{thm:imafcm} of piecewise $\mathcal C^{1+\alpha}$-diffeomorphic images of limit sets of sGDS. Let $\nu$ denote the $\mdim$-conformal measure associated with $R$. Then direct consequences of Thms.~\ref{thm:curvatureresult} and \ref{thm:imafcm} are:
  The average Minkowski content of both $\sset$ and $F$ exist and they are related by
  \[
  \widetilde{\mathcal M}(F)=\widetilde{\mathcal M}(\sset)\cdot\sum_{v\in V}\int_{\sset\cap Y_v}\lvert g_v'\rvert^{\mdim}\textup{d}\nu.
  \]
  If $R$ is non-lattice, then the Minkowski contents of both $\sset$ and $F$ exist and coincide with the respective average Minkowski contents.  
\end{remark}
\begin{remark}\label{rmk:Minkmeasurablelattice}
  From the result of Prop.~\ref{thm:conformalMinkowski} we can explicitly construct $\mathcal C^{1+\alpha}$-diffeomorphisms which map limit sets of lattice sGDS to Minkowski measurable limit sets of lattice cGDS. In fact, for every limit set $\sset$ of a lattice sGDS $R$ there exist $\mathcal C^{1+\alpha}$-diffeomorphisms $g$ such that $g(\sset)$ is Minkowski measurable:
  Assume that $\sset\subseteq[0,1]$ and that the geometric potential function $\ze$ associated with $R$ is lattice. Let $\aaa>0$ be maximal such that the range of $\ze$ is contained in $\aaa\mathbb Z$. Let $\nu$ denote the $\mdim$-conformal measure associated with $R$. Define $\widetilde{g}\colon\mathbb R\to\mathbb R$, $\widetilde{g}(x)\defeq\nu((-\infty,x])$ to be the distribution function of $\nu$. For $n\in\mathbb N$ define the function $g_n\colon[-1,\infty)\to\mathbb R$ by
  \[
  g_n(x)\defeq\int_{-1}^x\left(\widetilde{g}(r)(\ee^{\mdim\aaa n}-1)+1\right)^{-1/\mdim}\textup{d}r
  \]
  and set $F^n\defeq g_n(K)$. Then for every $n\in\mathbb N$ we have $\underline{\mathcal{M}}(F^n)=\overline{\mathcal{M}}(F^n)$ and $\underline{C}_0^f(F^n,\mathbb R)=\overline{C}_0^f(F^n, \mathbb R)$.
The proof of this statement has been given in \cite[Cor.~2.18(iii)]{KesKom} for self-conformal sets. For limit sets of cGDS the proof follows through by using Prop.~\ref{thm:conformalMinkowski} and thus, we are not going to repeat it here.
Notably, the sets $F^n$ do not only provide examples of Minkowski measurable limit sets of lattice cGDS but they also provide examples of sets for which $\mathcal M(F^n)$ and $C_0^f(F^n)$ exist but $C_1^f(F,\cdot)$ and $C_0^f(F,\cdot)$ do not exist (see Thm.~\ref{thm:imafcm}).
\end{remark}

We end this section with addressing conjectures from \cite{Dundee}.
\begin{remark}[On two conjectures by M.\,L.\,Lapidus from 1993]\label{rmk:conjLapidus}~\\[-3ex]
  \begin{enumerate}
    \item\label{it:conjss} Conjecture 3 in \cite{Dundee} states that under the OSC a non-degenerate self-similar set in $\mathbb R^d$ is Min\-kowski measurable if and only if the associated sIFS is non-lattice. This conjecture was proven to be correct in space dimension $d=1$ in \cite{Falconer_Minkowski}, \cite{Dundee}, \cite{Lapidus_Frankenhuysen_Springer}, \cite{LapPom} under the assumption that the feasible open set is connected. 
      In higher dimensional ambient spaces it was proven in \cite{Gatzouras} that a self-similar set arising from a non-lattice sIFS is Minkowski measurable, without any further assumptions on the feasible open set. Thus, the results from \cite{Gatzouras} fully prove one direction of the conjectured equivalence. The other direction, namely that one has Minkowski non-measurability in the lattice situation, is still an open problem in higher dimensional ambient spaces.  

With Cor.~\ref{cor:sGDSMink} we have seen that the Minkowski content of a limit set of an sGDS in $\mathbb R$ exists if and only if the sGDS is non-lattice. Thus, Cor.~\ref{cor:sGDSMink} shows that \cite[Conj.~3]{Dundee} is also valid for the more general class of limit sets of sGDS in $\mathbb R$ and in this way provides an important extension to the result of \cite{Falconer_Minkowski}, \cite{Dundee}, \cite{Lapidus_Frankenhuysen_Springer}, \cite{LapPom}. Moreover, Cor.~\ref{cor:sGDSMink} also allows to consider self-similar systems where the OSC is satisfied with disconnected feasible open sets and hence provides a new result in the lattice situation for self-similar sets (see Sec.~\ref{sec:disconnected}).
    \item\label{rmk:conjLapidus:Conj4} In the same paper, \cite{Dundee}, a similar conjecture is posed for so-called `approximately' self-similar sets, namely \cite[Conj.~4]{Dundee}. 
      A precise definition of an `approximately' self-similar set is not given, however, limit sets of Fuchsian groups of Schottky type are mentioned as examples. These can be represented as limit sets of cGDS (see Sec.~\ref{sec:Fuchsian}). It is well known that such systems are always non-lattice (see e.\,g.\ \cite[Part II]{Lalley}). Combined with \cite[Cor.~2.3]{LapPom}, Equation \eqref{eq:Mcnonlattice} thus verifies \cite[Conj.~4]{Dundee} for limit sets of Fuchsian groups of Schottky type. This situation will be investigated in more detail in Sec.~\ref{sec:Fuchsian}.

	In this spirit we view limit sets of cGDS in general as being `approximately' self-similar since conformal maps locally behave like similarities. Its subclass of self-conformal sets has already been treated in \cite{KesKom}. The results of \cite{KesKom} combined with \cite[Cor.~2.3]{LapPom} provide a negative answer to \cite[Conj.~4]{Dundee} for such sets (see also \cite[Ex.~2.20]{KesKom}). Note that \cite[Thm.~2.12]{KesKom} combined with \cite[Cor.~2.3]{LapPom} in particular shows that there exist fractal strings with lattice self-conformal boundary for which the asymptotic second term of the eigenvalue counting function $N(\lambda)$ of the Laplacian (in the sense of \cite{LapPom}) is monotonic.
	As self-conformal sets are special types of limit sets of cGDS, the results from \cite{KesKom} already imply that Minkowski measurability of the limit set of a cGDS is not equivalent to the cGDS being non-lattice. However, \eqref{eq:Mcnonlattice} shows the validity of one implication, namely that limit sets of non-lattice cGDS are Minkowski measurable.
	\end{enumerate}
\end{remark}

\section{Examples of Limit Sets of cGDS}\label{sec:examples}

We now present classes of systems which can be represented by a cGDS and illustrate our results for such systems. We especially focus on sets which cannot be treated with the previously known results from the literature.

\subsection{cGDS derived from a cIFS}\label{sec:CIFS_incidence}
A cIFS $\Psi\defeq(\psi_1,\ldots,\psi_N)$ has got the property that every function $\psi_i$ can be concatenated with any other function $\psi_j$ for $i,j\in\{1,\ldots,N\}$. Here we define a cGDS in that we additionally put transition rules on $\Psi$. This is done by defining an $N\times N$ matrix $A'\defeq(A'_{i,j})_{i,j\in\{1,\ldots,N\}}$ with entries $0,1$ which determines which functions may follow a given function, i.\,e.\ $A'_{i,j}=1$ if and only if $\psi_i\circ\psi_j$ is allowed. The system $(\Psi,A')$ then gives rise to a cGDS by setting $V\defeq\{1,\ldots,N\}$, $E\defeq \{1,\ldots,M\}$, where $M\defeq\sum_{i,j=1}^N A'_{i,j}$ and where for all $v,v'\in V$ with $A'_{v,v'}=1$ there exists an edge $e\in E$ such that $i(e)=v$ and $t(e)=v'$.	
\begin{example}\label{ex:CIFSincidence}
  For $i\in\{1,2,3\}$ define $\psi_i\colon[0,1]\to[0,1]$ by setting $\psi_1(x)\defeq x/4$, $\psi_2(x)\defeq x/4 +3/8$ and $\psi_3(x)\defeq x/4+3/4$ and set
  \[
  A'\defeq \left(\arraycolsep=2.5pt
  \begin{array}{ccc}
    1 & 0 & 1\\
    0 & 0 & 1\\
    1 & 1 & 1
  \end{array}\right).
  \]
  A corresponding sGDS is given by $V\defeq\{1,2,3\}$, $E\defeq\{1,\ldots,6\}$,
  \[
  i(e)\defeq \begin{cases}
    1 & \hspace{-0.2cm}:e\in\{1,2\}\\
    2 & \hspace{-0.2cm}:e=3\\
    3 & \hspace{-0.2cm}:e\in\{4,5,6\},
  \end{cases}\ 
  t(e)\defeq \begin{cases}
    1 & \hspace{-0.2cm}:e\in\{1,4\}\\
    2 & \hspace{-0.2cm}:e=5\\
    3 & \hspace{-0.2cm}:e\in\{2,3,6\},
  \end{cases}\ 
  A\defeq\left(\arraycolsep=2.5pt
  \begin{array}{cccccc}
    1&1&0&0&0&0\\
    0&0&0&1&1&1\\
    0&0&0&1&1&1\\
    1&1&0&0&0&0\\
    0&0&1&0&0&0\\
    0&0&0&1&1&1\\
  \end{array}\right),
  \]
  $X_v\defeq\psi_v([0,1])$ for $v\in V$ and 
  \begin{align*}
    &\phi_1\colon X_1\stackrel{\psi_1}{\longrightarrow} X_1 \qquad
    \phi_3\colon X_3\stackrel{\psi_2}{\longrightarrow} X_2 \qquad					
    \phi_5\colon X_2\stackrel{\psi_3}{\longrightarrow} X_3\\
    &\phi_2\colon X_3\stackrel{\psi_1}{\longrightarrow} X_1 \qquad
    \phi_4\colon X_1\stackrel{\psi_3}{\longrightarrow} X_3 \qquad					
    \phi_6\colon X_3\stackrel{\psi_3}{\longrightarrow} X_3. 	
  \end{align*}
  Here, $r=1/4$. For determining the average (local) Minkowski content of the limit set $F$ of the sGDS, we apply Thm.~\ref{thm:similars} and \eqref{eq:Mcaverage} and thus need to find the primary gaps. Observe that
  \[\begin{array}{lll}
  \langle\bij[1]\rangle = [0,1/16],&
  \langle\bij[2]\rangle = [3/16,1/4],&
  \langle\bij[3]\rangle = [9/16,5/8],\\
  \langle\bij[4]\rangle = [3/4,13/16],&
  \langle\bij[5]\rangle = [57/64,29/32]\quad \text{and}&
  \langle\bij[6]\rangle = [15/16,1].
  \end{array}\]
  Thus,
  \[
  L^1=\underbrace{\left(\frac{1}{16},\frac{3}{16}\right)}_{\eqdef L^{1,1}},\quad L^2=\varnothing\quad \text{and}\quad L^3=\underbrace{\left(\frac{13}{16},\frac{57}{64}\right)}_{\eqdef L^{3,1}}\cup\underbrace{\left(\frac{29}{32},\frac{15}{16}\right)}_{\eqdef L^{3,2}}.	
  \]
  The primary gaps $L^{1,1},L^{3,1}$ and $L^{3,2}$ are illustrated in Fig.~\ref{fig:incidence}.
  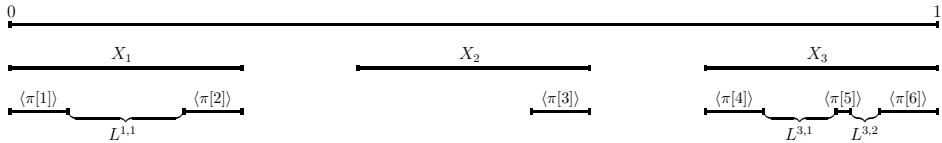
\begin{figure}[bb]
    \resizebox{\textwidth}{!}{
      \setlength{\unitlength}{0.348cm}
      \begin{picture}(65,10)
	\linethickness{0.4ex}
	\put(-0.2,9.5){\Large{0}}
	\put(63.7,9.5){\Large{1}}
	\put(0,9){\line(64,0){64}}
	\put(0,8.8){\line(0,1){0.4}}
	\put(64,8.8){\line(0,1){0.4}}
	
	\put(7,6.6){\Large{$X_1$}}
	\put(31,6.6){\Large{$X_2$}}
	\put(55,6.6){\Large{$X_3$}}
	\put(0,6){\line(64,0){16}}
	\put(0,5.8){\line(0,1){0.4}}
	\put(16,5.8){\line(0,1){0.4}}
	\put(24,6){\line(64,0){16}}
	\put(24,5.8){\line(0,1){0.4}}
	\put(40,5.8){\line(0,1){0.4}}
	\put(48,6){\line(64,0){16}}
	\put(48,5.8){\line(0,1){0.4}}
	\put(64,5.8){\line(0,1){0.4}}
	
	\put(0.65,3.6){\large{$\langle\bij[1]\rangle$}}
	\put(12.65,3.6){\large{$\langle\bij[2]\rangle$}}
	\put(36.65,3.6){\large{$\langle\bij[3]\rangle$}}
	\put(48.65,3.6){\large{$\langle\bij[4]\rangle$}}
	\put(56.15,3.6){\large{$\langle\bij[5]\rangle$}}
	\put(60.65,3.6){\large{$\langle\bij[6]\rangle$}}
	\put(0,3){\line(64,0){4}}
	\put(0,2.8){\line(0,1){0.4}}
	\put(4,2.8){\line(0,1){0.4}}
	\put(12,3){\line(64,0){4}}
	\put(12,2.8){\line(0,1){0.4}}
	\put(16,2.8){\line(0,1){0.4}}
	\put(36,3){\line(64,0){4}}
	\put(36,2.8){\line(0,1){0.4}}
	\put(40,2.8){\line(0,1){0.4}}
	\put(48,3){\line(64,0){4}}
	\put(48,2.8){\line(0,1){0.4}}
	\put(52,2.8){\line(0,1){0.4}}
	\put(57,3){\line(64,0){1}}
	\put(57,2.8){\line(0,1){0.4}}
	\put(58,2.8){\line(0,1){0.4}}
	\put(60,3){\line(64,0){4}}
	\put(60,2.8){\line(0,1){0.4}}
	\put(64,2.8){\line(0,1){0.4}}
	
	\put(4.1,2.9){$\underbrace{\hspace{2.7144cm}}$}
	\put(52.1,2.9){$\underbrace{\hspace{1.6704cm}}$}
	\put(58.1,2.9){$\underbrace{\hspace{0.6264cm}}$}
	\put(6.8,1){\Large{$L^{1,1}$}}
	\put(53.5,1){\Large{$L^{3,1}$}}
	\put(58,1){\Large{$L^{3,2}$}}
    \end{picture}}
    \caption{Primary gaps of the cGDS from Ex.~\ref{ex:CIFSincidence}.}
    \label{fig:incidence}
  \end{figure}
  Another quantity in the formula of Thm.~\ref{thm:similars} is the eigenfunction $\eigenf_{-\mdim\xi}$ of the Perron-Frobenius operator $\mathcal{L}_{-\mdim\xi}$ (see Sec.~\ref{sec:PF}), where $\mdim$ denotes the Min\-kowski dimension of $F$ and $\xi$ is the geometric potential function associated with $\Phi$. 
  In order to determine $\eigenf_{-\mdim\xi}$, we first determine the measure $\nu_{-\mdim\xi}$. This is done by solving the linear system of equations which arises by combining the following three facts:
  For $e\in E$ the defining equation for $\nu_{-\mdim\xi}$ implies that $\nu_{-\mdim\xi}([ee'])=4^{-\mdim}\cdot\nu_{-\mdim\xi}([e'])$ for every $e'\in T_{i(e)}$, 
  $\nu_{-\mdim\xi}([e])=\sum_{e'\in T_{i(e)}}\nu_{-\mdim\xi}([ee'])$ and 
  $\sum_{e\in E}\nu_{-\mdim\xi}([e])=1$. 
  The resulting measure $\nu_{-\mdim\xi}$ satisfies
  \begin{align*}
    \nu_{-\mdim\xi}([1])&=\nu_{-\mdim\xi}([4])=(3\cdot 4^{\mdim}-4^{-\mdim})^{-1},\\
    \nu_{-\mdim\xi}([2])&=\nu_{-\mdim\xi}([3])=\nu_{-\mdim\xi}([6])=(4^{\mdim}-1)\cdot\nu_{-\mdim\xi}([1]),\\
    \nu_{-\mdim\xi}([5])&=(1-4^{-\mdim})\cdot \nu_{-\mdim\xi}([1]).
  \end{align*}
  To determine $\eigenf_{-\mdim\xi}$, we use the approximation argument from \eqref{eq:convergenceperron}. We let $\mathbf 1$ denote the constant one-function on $E_A^{\infty}$. Since $\mathcal{L}_{-\mdim\xi}^n \mathbf{1}(\omneu)=\sum_{\om\in T_{v}^n}r_{\om}^{\mdim}$ for all $\omneu\in I_v^{\infty}$ and $v\in V$, it follows that $\eigenf_{-\mdim\xi}$ is constant on one-cylinders. Now combining the fact that the eigenvalue $\eigenv_{-\mdim\xi}$ is equal to one, that $\mathcal{L}_{-\mdim\xi}\eigenf_{-\mdim\xi}=\eigenv_{-\mdim\xi}\eigenf_{-\mdim\xi}$ and that $\int\eigenf_{-\mdim\xi}\textup{d}\nu_{-\mdim\xi}=1$, we obtain
  \begin{equation*}
    \begin{array}{rcll}
      \eigenf_{-\mdim\xi}(\om^1)&=&\displaystyle{\frac{3-4^{-2\mdim}}{-2\cdot 4^{-\mdim}+6-4^{\mdim}}}\qquad&\text{for}\ \om^1\in I_1^{\infty},\\[2ex]
      \eigenf_{-\mdim\xi}(\om^2)&=&(1-4^{-\mdim})\cdot\eigenf_{-\mdim\xi}(\om^1)&\text{for}\ \om^2\in I_2^{\infty}\\[1ex]
      \eigenf_{-\mdim\xi}(\om^3)&=&(4^{\mdim}-1)\cdot\eigenf_{-\mdim\xi}(\om^1)& \text{for}\ \om^3\in I_3^{\infty}.
    \end{array}
  \end{equation*}
  From the above evaluations we additionally infer that the Minkowski dimension $\mdim$ is the unique positive root of the function
  \[
  x\mapsto 4^{-x}-4^{-2x}+2-4^{x}.
  \]
  With $\entro(\mu_{-\mdim\xi})=\mdim\ln4$ we altogether obtain from Thm.~\ref{thm:similars} and \eqref{eq:Mcaverage} that 
  \[
  \widetilde{\mathcal{M}}(F)=\frac{2^{1-\mdim}\cdot(3-4^{-2\mdim})}{(1-\mdim)\mdim(6-2\cdot 4^{-\mdim}-4^{\mdim})\ln 4}
  \left(
  8^{-\mdim}+(4^{\mdim}-1)
  \left(\left(\frac{5}{64}\right)^{\mdim}+32^{-\mdim}\right)
  \right).
  \]
From Thm.~\ref{thm:curvatureresult} we conclude that 
	\[
	\widetilde{C}_1^f(F,\cdot)=\widetilde{\mathcal M}(F)\cdot \nu(\cdot)\quad\text{and}\quad
	\widetilde{C}_0^f(F,\cdot)=\frac{1-\mdim}{2}\widetilde{\mathcal M}(F)\cdot \nu(\cdot),
	\]
	where $\nu$ denotes the $\mdim$-conformal measure associated with $\Phi$.
  Since $\xi=\ln4\cdot\mathbf{1}$ is lattice, Cor.~\ref{cor:sGDSMink} moreover implies that the Minkowski content of $F$ does not exist.
\end{example}

\subsection{Conformal Iterated Function Systems with disconnected Feasible Open Set}\label{sec:disconnected}
By definition, a cIFS acting on $X$ needs to satisfy the OSC with $\inte(X)$ as a feasible open set. If we allow the OSC to be satisfied with a different feasible open set, then often the system can still be represented by a cGDS.
\begin{example}\label{ex:disconnected}
  For $i\in\{1,2,3\}$ define $\psi_i\colon[0,1]\to[0,1]$ by $\psi_1(x)\defeq x/3$, $\psi_2(x)\defeq x/3+2/3$ and $\psi_3(x)\defeq x/9+1/9$ and set $\Psi\defeq(\psi_1,\psi_2,\psi_3)$. Then $\Psi$ is not a cIFS in our sense since the open set condition is not satisfied with $(0,1)$ as the feasible open set. (Even though the OSC is satisfied for $(0,1/3)\cup(2/3,1)$.) However, $\Psi$ can be represented by an sGDS as follows.
  Set $V\defeq\{1,2\}$, $E\defeq\{1,\ldots,6\}$, 
  \[
  i(e)\defeq\begin{cases}
  1 & \hspace{-0.2cm}:e\in\{1,\ldots,4\}\ \\
  2 & \hspace{-0.2cm}:e\in\{5,6\},
  \end{cases}\quad
  t(e)\defeq\begin{cases}
  1 & \hspace{-0.2cm}:e\in\{1,3,5\}\\
  2 & \hspace{-0.2cm}:e\in\{2,4,6\},\ 
  \end{cases}\quad 
  A\defeq \left(\arraycolsep=2.5pt
  \begin{array}{cccccc}
    1&1&1&1&0&0\\
    0&0&0&0&1&1\\
    1&1&1&1&0&0\\
    0&0&0&0&1&1\\
    1&1&1&1&0&0\\
    0&0&0&0&1&1\\
  \end{array}\right),
  \]
  $X_v\defeq \psi_v([0,1])$ for $v\in\{1,2\}$ and
  \[\begin{array}{lll}
    \phi_1\colon X_1\stackrel{\psi_1}{\longrightarrow} X_1,\qquad &\phi_3\colon X_1\stackrel{\psi_3}{\longrightarrow} X_1,\qquad &\phi_5\colon X_1\stackrel{\psi_2}{\longrightarrow} X_2,\\
    \phi_2\colon X_2\stackrel{\psi_1}{\longrightarrow} X_1, &\phi_4\colon X_2\stackrel{\psi_3}{\longrightarrow} X_1\ \ \text{and}\ &\phi_6\colon X_2\stackrel{\psi_2}{\longrightarrow} X_2.
  \end{array}\]%
  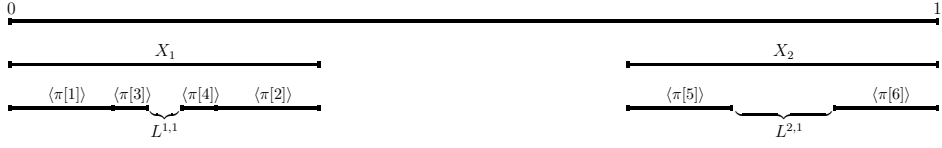
\begin{figure}[t]
    \resizebox{\textwidth}{!}{
      \setlength{\unitlength}{0.348cm}
      \begin{picture}(65,12)
	\linethickness{0.4ex}
	\put(-0.2,9.5){\Large{0}}
	\put(63.7,9.5){\Large{1}}
	\put(0,9){\line(64,0){64}}
	\put(0,8.8){\line(0,1){0.4}}
	\put(64,8.8){\line(0,1){0.4}}
	
	\put(10,6.6){\Large{$X_1$}}
	\put(52.6667,6.6){\Large{$X_2$}}
	\put(0,6){\line(64,0){21.3333}}
	\put(0,5.8){\line(0,1){0.4}}
	\put(21.3333,5.8){\line(0,1){0.4}}
	\put(42.6667,6){\line(64,0){21.3333}}
	\put(42.6667,5.8){\line(0,1){0.4}}
	\put(64,5.8){\line(0,1){0.4}}
	
	\put(2.6,3.6){\large{$\langle\bij[1]\rangle$}}
	\put(16.8222,3.6){\large{$\langle\bij[2]\rangle$}}
	\put(7.1111,3.6){\large{$\langle\bij[3]\rangle$}}
	\put(11.8519,3.6){\large{$\langle\bij[4]\rangle$}}
	\put(45.2667,3.6){\large{$\langle\bij[5]\rangle$}}
	\put(59.48889,3.6){\large{$\langle\bij[6]\rangle$}}
	\put(0,3){\line(64,0){7.1111}}
	\put(0,2.8){\line(0,1){0.4}}
	\put(7.1111,2.8){\line(0,1){0.4}}
	\put(7.1111,3){\line(64,0){2.3704}}
	\put(9.4815,2.8){\line(0,1){0.4}}
	\put(11.8519,3){\line(64,0){2.3704}}
	\put(11.8519,2.8){\line(0,1){0.4}}
	\put(14.2222,3){\line(64,0){7.1111}}
	\put(14.2222,2.8){\line(0,1){0.4}}
	\put(21.3333,2.8){\line(0,1){0.4}}
	\put(42.6667,3){\line(64,0){7.1111}}
	\put(42.6667,2.8){\line(0,1){0.4}}
	\put(49.7778,2.8){\line(0,1){0.4}}
	\put(56.8889,3){\line(64,0){7.1111}}
	\put(56.8889,2.8){\line(0,1){0.4}}
	\put(64,2.8){\line(0,1){0.4}}
	
	\put(9.6815,2.9){$\underbrace{\hspace{0.68cm}}$}
	\put(50.0778,2.9){$\underbrace{\hspace{2.3022cm}}$}
	\put(9.7,1){\Large{$L^{1,1}$}}
	\put(52.8,1){\Large{$L^{2,1}$}}
    \end{picture}}
    \caption{Primary gaps of the limit set of the cGDS from Ex.~\ref{ex:disconnected}.}
    \label{fig:disconnected}
  \end{figure}%
  Here, $r=1/3$, $L^{1,1}=(4/27,5/27)$ and $L^{2,1}=(7/9,8/9)$. See Fig.~\ref{fig:disconnected} for an illustration for this example. 
  That the eigenfunction $\eigenf_{-\mdim\xi}$ of the Perron-Frobenius operator $\mathcal{L}_{-\mdim\xi}$ with eigenvalue 1 is equal to the constant one-function $\mathbf 1$ on $E_A^{\infty}$ can be seen as follows.
  Firstly, $\mathcal{L}_{-\mdim\xi}\mathbf{1}= 2/3^{\mdim}+1/9^{\mdim}$ and secondly, $1=2/3^{\mdim}+1/9^{\mdim}$ which can be concluded from the fact that $0=P(-\mdim\xi)$, where $P$ denotes the topological pressure function (see \eqref{eq:pressuredef}). Thus, by Thm.~\ref{thm:similars} and \eqref{eq:Mcaverage}, we have
  \[
  \widetilde{\mathcal M}(F)=
  \frac{2^{1-\mdim}\cdot(27^{-\mdim}+9^{-\mdim})}{(1-\mdim)\entro(\mu_{-\mdim\xi})}.
  \]
As in the previous example $\widetilde{C}_1^f(F,\cdot)$ and $\widetilde{C}_0^f(F,\cdot)$ can be determined from the above equation by using Thm.~\ref{thm:curvatureresult}.
  Cor.~\ref{cor:sGDSMink} implies that the Minkowski content of $F$ does not exist, since the range of $\xi$ is contained in $\ln 3\cdot\mathbb Z$.\\
  Alternatively, one can determine the average Minkowski content of this example by using the results of \cite{Gatzouras}.
  However, if $\psi_1,\psi_2$ and $\psi_3$ were non-linear but conformal, then Thms.~\ref{thm:curvatureresult} and \ref{thm:conformalMinkowski} could be applied, whereas this case is not covered in \cite{Gatzouras}. 
\end{example}
		
\subsection{Markov Interval Maps}
For closed intervals $X_1,\ldots,X_N$ in $[0,1]$ with disjoint interior, $N\geq 2$, and $X:=\bigcup_{i=1}^N X_{i}$ we call a map $g\colon X\to[0,1]$ a \emph{Markov interval map} if 
\begin{enumerate}
\item $g|_{X_{i}}$ is expanding and there exists a $C^{1+\alpha}$-continuation to a neighbourhood of $X_{i}$ and
\item if $g(X_{i})\cap X_{j}\neq\varnothing$ then $X_{j}\subset g(X_{i})$ for $i,j\in\{1,\ldots,N\}$.
\end{enumerate}
For a representation by a cGDS, set $V\defeq\{1,\ldots,N\}$ and for $v\in V$ define $G_v\defeq\{v'\in V\mid X_{v'}\subseteq g(X_v)\}$. For every pair $(v,v')$, where $v\in V$ and $v'\in G_v$ introduce an edge $e=e(v,v')$ with $i(e)=v$ and $t(e)=v'$. Set $E\defeq\{e(v,v')\mid v\in V, v'\in G_v\}$ and define $\phi_e\colon X_{t(e)}\to X_{i(e)}$ by $\phi_e\defeq \left(g\vert_{X_{i(e)}}\right)^{-1}\vert_{X_{t(e)}}$ for $e\in E$.
Then the repeller of the Markov interval map coincides with the limit set of the corresponding cGDS.

\begin{example}\label{ex:Markovinterval}
  Set $X_1\defeq[0,1/4]$, $X_2\defeq[1/4,1/2]$, $X_3\defeq[2/3,1]$ and let the Markov interval map $g\colon \bigcup_{i=1}^3 X_i\to[0,1]$ be given by $g\vert_{X_1}(x)\defeq 5x/2$, $g\vert_{X_2}(x)\defeq 3x-1/2$ and $g\vert_{X_3}(x)\defeq 3x-2$. The graph of the Markov interval map $g$ is presented in Fig.~\ref{fig:Markovinterval_map}.
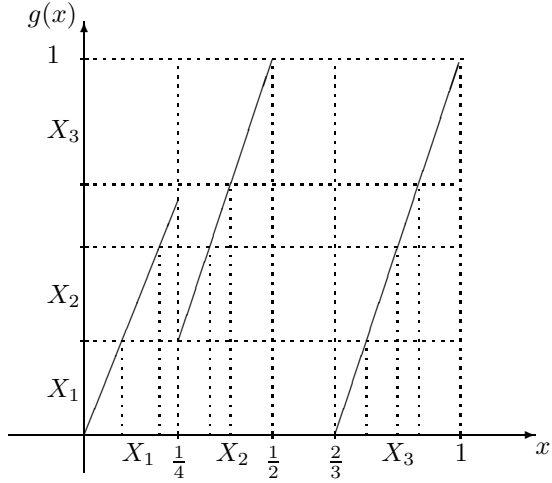
\begin{figure}[h!]
  \setlength{\unitlength}{5cm}
  \begin{picture}(90,1.45)
    
    \put(0.8,0.1){\vector(1,0){1.4}}
    \put(2.2,0.05){$x$}
    \put(1,0){\vector(0,1){1.2}}
    \put(0.85,1.2){$g(x)$}
    
    \put(1.23,0.01){$\frac{1}{4}$}
    \put(1.48,0.01){$\frac12$}
    \put(1.6467,0.01){$\frac23$}
    \put(1.985,0.03){$1$}
    \put(1.1,0.03){$X_1$}
    \put(1.35,0.03){$X_2$}
    \put(1.79,0.03){$X_3$}
    \put(1.25,0.09){\line(0,1){0.02}}
    \put(1.5,0.09){\line(0,1){0.02}}
    \put(1.6667,0.09){\line(0,1){0.02}}
    \put(2,0.09){\line(0,1){0.02}}
    
    \put(0.9,1.09){$1$}
    \put(0.9,0.2){$X_1$}
    \put(0.9,0.45){$X_2$}
    \put(0.9,0.89){$X_3$}
    \put(0.99,0.35){\line(1,0){0.02}}
    \put(0.99,0.6){\line(1,0){0.02}}
    \put(0.99,0.7667){\line(1,0){0.02}}
    \put(0.99,1.1){\line(1,0){0.02}}
    
    \put(1,0.1){\line(2,5){0.25}}
    \put(1.25,0.35){\line(1,3){0.25}}
    \put(1.6667,0.1){\line(1,3){0.33}}
    
    \multiput(1.25,0.1)(0,0.025){40}{\line(0,1){0.006}}
    \multiput(1.5,0.1)(0,0.025){40}{\line(0,1){0.006}}
    \multiput(1.6667,0.1)(0,0.025){40}{\line(0,1){0.006}}
    \multiput(2,0.1)(0,0.025){40}{\line(0,1){0.006}}
    
    \multiput(1.1,0.1)(0,0.025){10}{\line(0,1){0.004}}	
    \multiput(1.2,0.1)(0,0.025){20}{\line(0,1){0.004}}
    \multiput(1.3333,0.1)(0,0.025){20}{\line(0,1){0.004}}
    \multiput(1.3889,0.1)(0,0.025){27}{\line(0,1){0.004}}
    \multiput(1.75,0.1)(0,0.025){10}{\line(0,1){0.004}}
    \multiput(1.8333,0.1)(0,0.025){20}{\line(0,1){0.004}}
    \multiput(1.8889,0.1)(0,0.025){27}{\line(0,1){0.004}}
    
    \multiput(1,0.35)(0.025,0){40}{\line(1,0){0.007}}
    \multiput(1,0.6)(0.025,0){40}{\line(1,0){0.007}}
    \multiput(1,0.7667)(0.025,0){40}{\line(1,0){0.007}}
    \multiput(1,1.1)(0.025,0){41}{\line(1,0){0.007}}
  \end{picture}
  \caption{Graph of the Markov interval map from Ex.~\ref{ex:Markovinterval}.}
  \label{fig:Markovinterval_map}
\end{figure}
				
A corresponding sGDS is given by $V\defeq\{1,2,3\}$, $E\defeq\{1,\ldots,7\}$,
\[
i(e)\defeq\begin{cases}
1 &\hspace{-0.2cm}:e\in\{1,2\}\\
2 &\hspace{-0.2cm}:e\in\{3,4\}\\
3 &\hspace{-0.2cm}:e\in\{5,6,7\},\ 
\end{cases}\quad
t(e)\defeq\begin{cases}
1 & \hspace{-0.2cm}:e\in\{1,5\}\\
2 & \hspace{-0.2cm}:e\in\{2,3,6\}\ \\
3 & \hspace{-0.2cm}:e\in\{4,7\},
\end{cases}\quad
A\defeq\left(\arraycolsep=2.5pt
\begin{array}{ccccccc}
  1&1&0&0&0&0&0\\
  0&0&1&1&0&0&0\\
  0&0&1&1&0&0&0\\
  0&0&0&0&1&1&1\\
  1&1&0&0&0&0&0\\
  0&0&1&1&0&0&0\\
  0&0&0&0&1&1&1
\end{array}
\right),
\]
\begin{align*}
  &\phi_1\colon X_1 \stackrel{\left(g\vert_{X_1}\right)^{-1}}{\xrightarrow{\hspace*{0.9cm}}} X_1,\qquad
  \phi_2\colon X_2 \stackrel{\left(g\vert_{X_1}\right)^{-1}}{\xrightarrow{\hspace*{0.9cm}}} X_1,\\
  &\phi_3\colon X_2 \stackrel{\left(g\vert_{X_2}\right)^{-1}}{\xrightarrow{\hspace*{0.9cm}}} X_2,\qquad 
  \phi_4\colon X_3 \stackrel{\left(g\vert_{X_2}\right)^{-1}}{\xrightarrow{\hspace*{0.9cm}}} X_2,\\
  &\phi_5\colon X_1 \stackrel{\left(g\vert_{X_3}\right)^{-1}}{\xrightarrow{\hspace*{0.9cm}}} X_3,\qquad  
  \phi_6\colon X_2 \stackrel{\left(g\vert_{X_3}\right)^{-1}}{\xrightarrow{\hspace*{0.9cm}}} X_3,\qquad
  \phi_7\colon X_3 \stackrel{\left(g\vert_{X_3}\right)^{-1}}{\xrightarrow{\hspace*{0.9cm}}} X_3.
\end{align*}
Here, $r=3/4$.
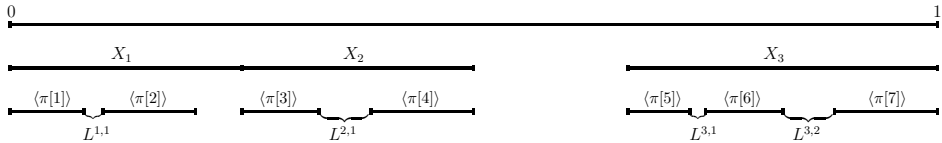
\begin{figure}[b]
  \resizebox{\textwidth}{!}{
    \setlength{\unitlength}{0.348cm}
    \begin{picture}(65,12)
      \linethickness{0.4ex}
      
      \put(-0.2,9.5){\Large{0}}
      \put(63.7,9.5){\Large{1}}
      \put(0,9){\line(64,0){64}}
      \put(0,8.8){\line(0,1){0.4}}
      \put(64,8.8){\line(0,1){0.4}}
      
      \put(7,6.6){\Large{$X_1$}}
      \put(23,6.6){\Large{$X_2$}}
      \put(52,6.6){\Large{$X_3$}}
      \put(0,6){\line(64,0){16}}
      \put(0,5.8){\line(0,1){0.4}}
      \put(16,5.8){\line(0,1){0.4}}
      \put(16,6){\line(64,0){16}}
      \put(32,5.8){\line(0,1){0.4}}
      \put(42.6667,6){\line(64,0){21.3333}}
      \put(42.6667,5.8){\line(0,1){0.4}}
      \put(64,5.8){\line(0,1){0.4}}
      
      \put(1.6,3.6){\large{$\langle\bij[1]\rangle$}}
      \put(8.2,3.6){\large{$\langle\bij[2]\rangle$}}
      \put(17.2,3.6){\large{$\langle\bij[3]\rangle$}}
      \put(27,3.6){\large{$\langle\bij[4]\rangle$}}
      \put(43.7,3.6){\large{$\langle\bij[5]\rangle$}}
      \put(49.2,3.6){\large{$\langle\bij[6]\rangle$}}
      \put(59.2,3.6){\large{$\langle\bij[7]\rangle$}}
      \put(0,3){\line(64,0){5.12}}
      \put(0,2.8){\line(0,1){0.4}}
      \put(5.12,2.8){\line(0,1){0.4}}
      \put(6.4,3){\line(64,0){6.4}}
      \put(6.4,2.8){\line(0,1){0.4}}
      \put(12.8,2.8){\line(0,1){0.4}}
      \put(16,3){\line(64,0){5.3333}}
      \put(16,2.8){\line(0,1){0.4}}
      \put(21.333,2.8){\line(0,1){0.4}}
      \put(24.8889,3){\line(64,0){7.1111}}
      \put(24.8889,2.8){\line(0,1){0.4}}
      \put(32,2.8){\line(0,1){0.4}}
      \put(42.6667,3){\line(64,0){4.2667}}
      \put(42.6667,2.8){\line(0,1){0.4}}
      \put(46.9337,2.8){\line(0,1){0.4}}
      \put(48,3){\line(64,0){5.3333}}
      \put(48,2.8){\line(0,1){0.4}}
      \put(53.3333,2.8){\line(0,1){0.4}}
      \put(56.8889,3){\line(64,0){7.1111}}
      \put(56.8889,2.8){\line(0,1){0.4}}
      \put(64,2.8){\line(0,1){0.4}}
      
      \put(5.16,2.9){\tiny{$\underbrace{\hspace{0.1cm}}$}}
      \put(21.4333,2.9){$\underbrace{\hspace{1.15cm}}$}
      \put(46.9,2.9){\tiny{$\underbrace{\hspace{0.1cm}}$}}
      \put(53.4333,2.9){$\underbrace{\hspace{1.15cm}}$}
      
      \put(5,1){\Large{$L^{1,1}$}}
      \put(22,1){\Large{$L^{2,1}$}}
      \put(46.9,1){\Large{$L^{3,1}$}}
      \put(54,1){\Large{$L^{3,2}$}}
      
  \end{picture}}
  \caption{Primary gaps for the limit set of Ex.~\ref{ex:Markovinterval}.}
  \label{fig:Markovgaps}
\end{figure}
For this example, we limit ourselves to determining and illustrating the primary gaps, since presenting the complete calculations would not provide any further insights.
The convex hulls of the projections of the cylinder sets are given by
\begin{equation*}
  \begin{array}{lll}
    \langle\bij[1]\rangle=[0,2/25],\qquad&  \langle\bij[2]\rangle=[1/10,1/5],\qquad&\\
    \langle\bij[3]\rangle=[1/4,1/3],& \langle\bij[4]\rangle=[7/18,1/2],&\\
    \langle\bij[5]\rangle=[2/3,11/15],\qquad& \langle\bij[6]\rangle=[3/4,5/6],&\langle\bij[7]\rangle=[8/9,1].
  \end{array}
\end{equation*}
Thus, the primary gaps are 
\[\begin{array}{ll}
L^{1,1}=(2/25,1/10),\qquad &L^{2,1}=(1/3,7/18),\\
L^{3,1}=(11/15,3/4)\ \ \text{and}\  &L^{3,2}=(5/6,8/9). 
\end{array}\]
They are illustrated in Fig.~\ref{fig:Markovgaps}.
This cGDS indeed is a non-lattice cGDS and hence $C_1^f(F,\cdot)$,  $C_0^f(F,\cdot)$ and the Minkowski content of $F$ exist by Thm.~\ref{thm:similars} and \eqref{eq:Mcnonlattice}.
\end{example}

\subsection{Lattice cGDS whose Limit Set is Minkowski Measurable}
An example of a lattice self-conformal set which is Minkowski measurable is given in \cite[Ex.~2.20]{KesKom}. 
In the following, we present an example of a Minkowski measurable limit set of a lattice cGDS which cannot be obtained via a cIFS. 
This adds to the observations concerning \cite[Conj.~4]{Dundee} that we discussed in Rem.~\ref{rmk:conjLapidus}\ref{rmk:conjLapidus:Conj4}.
To be more precise, the following example in conjunction with \cite[Cor.~2.3]{LapPom} shows the existence of fractal strings, that have a limit set of a lattice cGDS for boundary, for which the asymptotic second term of the eigenvalue counting function $N(\lambda)$ of the Laplacian (in the sense of \cite{LapPom}) is monotonic. This shows that the statement of \cite[Conj.~4]{Dundee} is not valid for limit sets of cGDS in $\mathbb R$.

\begin{example}\label{ex:thexample}
	Let $K\subseteq[0,1]$ denote the limit set of the sGDS given in Ex.~\ref{ex:CIFSincidence}. Let $\mdim$ denote its Minkowski dimension and let $\nu$ denote the associated $\mdim$-conformal measure. Let $\widetilde{g}\colon\mathbb R\to\mathbb R$ denote the distribution function of $\nu$, i.\,e.\ $\widetilde{g}(x)\defeq\nu((-\infty,x])$ for $x\in\mathbb R$. For $n\in\mathbb N$ define the function $g_n\colon[-1,\infty)\to\mathbb R$ by
	\[
		g_n(x)\defeq\int_{-1}^{x}\left(\widetilde{g}(r)(3^{n\mdim}-1)+1\right)^{-1/\mdim}\textup{d}r
	\]
	and set $F^n\defeq g_n(K)$. Then we have $\underline{\mathcal M}(F^n)=\overline{\mathcal M}(F^n)$, although $\underline{\mathcal M}(K)<\overline{\mathcal M}(K)$. This is a consequence of Cor.~\ref{cor:sGDSMink} and Rem.~\ref{rmk:Minkmeasurablelattice}.
\end{example}

\subsection{Limit Sets of Fuchsian Groups of Schottky Type}\label{sec:Fuchsian}
Here, we give a very brief introduction to limit sets of Fuchsian groups of Schottky type. For background and proofs of the statements below, we refer the reader to \cite{Beardon}, \cite{Nicholls}.

We let $\mathbb H\defeq\{z\in\mathbb C\mid \Im(z)>0\}$ denote the upper half plane in $\mathbb C$, where $\Im(z)$ denotes the imaginary part of $z\in\mathbb C$. We fix $n\in\mathbb N$ with $n\geq 2$ and set $V\defeq\{\pm 1,\ldots,\pm n\}$. We let $(B_v)_{v\in V}$ denote a family of pairwise disjoint closed Euclidean unit balls in $\mathbb C$ intersecting the real line $\mathbb R$ orthogonally and let $g_v$ denote the unique hyperbolic conformal orientation preserving automorphism of $\mathbb H$ which maps the side $s_{-v}\defeq\mathbb H\cap\partial B_{-v}$ to the side $s_v\defeq\mathbb H\cap\partial B_v$. (Note that $g_v$ is a M\"obius transformation which is obtained on concatenating the inversion at the circle $\partial B_{-v}$ with the reflection at the line $\Re(z)=d_v$, where $d_v=(c_v+c_{-v})/2$ is the midpoint of the line segment joining the centres $c_{-v}$ and $c_v$ of the balls $B_{-v}$ and $B_v$ and $\Re(z)$ denotes the real part of $z\in \mathbb C$.)
Then $\{g_v\mid v\in V\}$ is a symmetric set of generators of the Fuchsian group $G\defeq\langle\{g_v\mid v\in V\}\rangle$ and $G$ will be referred to as a \emph{Fuchsian group of Schottky type}. Associated to $G$ is a limit set $L(G)\subset\mathbb R\cap\bigcup_{v\in V}B_v$ which is defined to be the set of all accumulation points (with respect to the Euclidean metric on $\overline{\mathbb H}\defeq\mathbb H\cup\mathbb R\cup\{\infty\}$) of the $G$-orbit $G(z)\defeq\{g(z)\mid g\in G\}$ for an arbitrary $z\in\mathbb H$. 

Such a limit set can be represented as a limit set of a cGDS in the following way:
For defining the directed multigraph we set the set of vertices to be $V$, define $E\defeq\{(v,v')\in V^2\mid v'\neq -v\}$ to be the set of edges, $t((v,v'))\defeq v$ and $i((v,v'))\defeq v'$.
The incidence matrix $A$ is given by $A_{e,e'}=1$ if $t(e)=i(e')$ and $A_{e,e'}=0$ else. It is aperiodic and irreducible, which can be seen as follows.
Let $e,e'''\in E$ denote two arbitrary edges. The condition that $n\geq 2$ implies that there exist at least two vertices $v\in V\setminus\{-t(e),-i(e''')\}$. Fix $v$ as such. Since $v\neq -t(e)$ there exists an edge $e'\in E$ with $i(e')=t(e)$ and $t(e')=v$ and likewise, there exists an edge $e''\in E$ with $i(e'')=v$ and $t(e'')=i(e''')$. Thus, $A^{(3)}_{e,e'''}>0$. 
For $v\in V$ we set $X_v\defeq B_v\cap\mathbb R$ and note that the maps $g_v$ can be continuously extended to $\overline{\mathbb H}$. We denote this extension also by $g_v$. For each $e=(t(e),i(e))\in E$ we set
\[
	\phi_e\colon X_{t(e)} \xrightarrow{g_{i(e)}} X_{i(e)}.
\]
Since each $g_v$ is a M\"obius transformation with singularity in $X_{-v}$, the map $\phi_e$ extends to an analytic $\mathcal C^{1+\alpha}$-diffeomorphism on an open connected neighbourhood $W_{t(e)}$ of $X_{t(e)}$, for some $\alpha\in(0,1]$. Moreover, the maps $\phi_e$ are strict contractions by construction. 
That the limit set $L(G)$ of the Fuchsian group coincides with the limit set of the above constructed cGDS is shown in \cite[Thm.~5.1.6]{Urbanski_Buch}.
By \cite[Part II]{Lalley} the associated geometric potential function is non-lattice. Therefore, we obtain the following corollary from Thm.~\ref{thm:curvatureresult}:

\begin{corollary}\label{cor:Fuchsian}
	The local Minkowski content and the local fractal Euler characteristic of a limit set of a Fuchsian group of Schottky type always exist. In particular, a limit set of a Fuchsian group of Schottky type is always Minkowski measurable.
\end{corollary}

Note that the above corollary proves \cite[Conj.~4]{Dundee} for limit sets of Fuchsian groups of Schottky type.
\begin{example}\label{ex:Fuchsian}
	In this example we want to show how a typical limit set of a Fuchsian group of Schottky type can be represented as a cGDS. We set $V\defeq\{\pm 1,\pm 2\}$ and define $B_{-2},B_{-1},B_{1}$ and $B_{2}$ to be the closed unit balls with respective centres $-5, -2, 2$ and $5$. Then the maps $g_v\colon \mathbb H\to\mathbb H$ are given by
	\begin{align*}
		g_{-2}(z)=\frac{-5z+24}{z-5},\ 
		g_{-1}(z)=\frac{-2z+3}{z-2},\ 
		g_{1}(z)=\frac{2z+3}{z+2}\ \text{and}\ 
		g_{2}(z)=\frac{5z+24}{z+5}
	\end{align*}
and $G\defeq\langle\{g_v\mid v\in V\}\rangle$ is the Fuchsian group of Schottky type. For a representation by a cGDS we set
\[
	X_{-2}\defeq[-6,-4],\quad X_{-1}\defeq[-3,-1],\quad X_1\defeq[1,3]\quad\text{and}\quad X_2\defeq [4,6].
\]
The set of edges is given by $E\defeq\{(v,v')\in V^2\mid v'\neq -v\}$, $t((v,v'))=v$, $i((v,v'))=v'$ and the family of maps $\phi_e$ for $e\in E$ is given by
\[\begin{array}{rrrrrr}
	\phi_{(-2,-2)}:&\hspace{-0.25cm}X_{-2}\xrightarrow{g_{-2}} X_{-2},\ \
	&\phi_{(-2,-1)}:&\hspace{-0.25cm}X_{-2}\xrightarrow{g_{-1}} X_{-1},\ \  
	&\phi_{(-2,1)}:&\hspace{-0.25cm}X_{-2}\xrightarrow{g_{1}} X_{1},\\
	\phi_{(-1,-2)}:&\hspace{-0.25cm}X_{-1}\xrightarrow{g_{-2}} X_{-2},\ \   
	&\phi_{(-1,-1)}:&\hspace{-0.25cm}X_{-1}\xrightarrow{g_{-1}} X_{-1},\ \   
	&\phi_{(-1,2)}:&\hspace{-0.25cm}X_{-1}\xrightarrow{g_{2}} X_{2},\\
	\phi_{(1,-2)}:&\hspace{-0.25cm}X_{1}\xrightarrow{g_{-2}} X_{-2},\ \ 
	&\phi_{(1,1)}:&\hspace{-0.25cm}X_{1}\xrightarrow{g_{1\hphantom{-}}} X_{1\hphantom{-}},\ \ 
	&\phi_{(1,2)}:&\hspace{-0.25cm}X_{1}\xrightarrow{g_{2}} X_{2},\\
	\phi_{(2,-1)}:&\hspace{-0.25cm}X_{2}\xrightarrow{g_{-1}} X_{-1},\ \ 
	&\phi_{(2,1)}:&\hspace{-0.25cm}X_{2}\xrightarrow{g_{1\hphantom{-}}} X_{1\hphantom{-}}, \ \
	&\phi_{(2,2)}:&\hspace{-0.25cm}X_{2}\xrightarrow{g_{2}} X_{2}.
\end{array}\]
The incidence matrix $A$ is a $12\times 12$ matrix which contains exactly three ones in every row and every column.
\end{example}

\section{Preliminaries}\label{sec:preliminaries}

We now provide some background information and auxiliary results for proving our main theorems, which we presented in Sec.~\ref{sec:mainresults}.

\subsection{Perron-Frobenius Theory and the Geometric Potential Function}\label{sec:PF}
In order to provide the necessary background to define the constants in our main statements and also to set up the tools needed in the proofs we now recall some facts from the Perron-Frobenius theory. 
For this, we are going to make use of results from \cite{Urbanski_Buch} which were obtained for conformal graph directed Markov systems (cGDMS), see Rem.~\ref{rmk:GDMS}, which are finitely primitive. A cGDMS is called \emph{finitely primitive}, if there exists an $n\in\mathbb N$ such that for all $e,e'\in E$ there exists an $\om\in E_A^n$ for which $e\om e'\in E_A^*$.

\begin{remark}\label{rmk:finitelyprimitive}
	A cGDS with aperiodic and irreducible incidence matrix is a finitely primitive cGDMS.
\end{remark}

In this subsection we always assume that the incidence matrix $A$ is aperiodic and irreducible.
Recall from Sec.~\ref{sec:Notation} that we equip $E_A^{\infty}$ as defined in \eqref{eq:EAinfty} with the sub-topology of the product topology of the discrete topologies of $E$ and let $\mathcal C(E_A^{\infty})$ denote the set of real-valued continuous functions on $E_A^{\infty}$.
For $f\in\mathcal C(E_A^{\infty})$, $\alpha\in(0,1)$ and $n\in\mathbb N\cup\{0\}$ we define 
\begin{align*}
  \text{var}_n(f)&\defeq\sup\{\lvert f(\om)-f(\omneu)\rvert\mid \om, \omneu\in E_A^{\infty}\ \text{and}\ \om_i=\omneu_i\ \text{for}\ i\in\{1,\ldots,n\}\},\\
  \lvert f\rvert_{\alpha}&\defeq\sup_{n\geq 0}\frac{\text{var}_n(f)}{\alpha^{n}}\ \qquad\text{and}\\
  \mathcal F_{\alpha}( E_A^{\infty})&\defeq\{f\in\mathcal C( E_A^{\infty})\mid \lvert f\rvert_{\alpha}<\infty\}.
\end{align*}
Elements of $\mathcal F_{\alpha}( E_A^{\infty})$ are called \emph{$\alpha$-H\"older continuous} functions on $ E_A^{\infty}$.
The space $\mathcal F_{\alpha}( E_A^{\infty})$ endowed with the norm $\|\cdot\|_{\alpha}\defeq |\cdot|_{\alpha}+\|\cdot\|$, where $\|\cdot\|$ denotes the supremum norm, is a Banach space.

\begin{remark}\label{rmk:alpha}
	The geometric potential function $\xi$ associated with a cGDS $\Phi\defeq\{\phi_e\}_{e\in E}$ satisfies $\xi\in\mathcal F_{\widetilde{\alpha}}(E_A^{\infty})$ for some $\widetilde{\alpha}\in(0,1)$. To see this, we let $r\in(0,1)$ denote a common Lipschitz constant of $\phi_e$ for $e\in E$. Because of the $\alpha$-H\"older continuity of $\phi'_e$, we obtain that there exists a constant $c\in\mathbb R$ such that for every $n\in\mathbb N$ we have $\text{var}_n(\xi)\leq cr^{\alpha(n-1)}$ and $\text{var}_0(\xi)<\infty$. Thus, $\xi\in\mathcal F_{\widetilde{\alpha}}( E_A^{\infty})$, where $\widetilde{\alpha}\defeq r^{\alpha}\in(0,1)$.
\end{remark}
For $f\in\mathcal{C}( E_A^{\infty})$ define the \emph{Perron-Frobenius operator} $\mathcal L_f\colon\mathcal C( E_A^{\infty})\to\mathcal C( E_A^{\infty})$ by
\begin{equation}\label{PerronFrobenius}
	\mathcal L_f \psi(\om)\defeq\sum_{\omneu:\sigma \omneu=\om} \ee^{f(\omneu)}\psi(\omneu)
\end{equation}
for $\om\in E_A^{\infty}$ and let $\mathcal{L}_{ f}^*$ be the dual of $\mathcal{L}_{ f}$ acting on the set of Borel probability measures on $ E_A^{\infty}$.
By \cite[Thm.~2.16 and Cor.~2.17]{Walters_convergence} and \cite[Thm.~1.7]{Bowen_equilibrium}, for each real-valued H\"older continuous $f\in\mathcal F_{\alpha}( E_A^{\infty})$, some $\alpha\in(0,1)$, there exists a unique Borel probability measure $\nu_{ f}$ on $ E_A^{\infty}$ such that $\mathcal L_{ f}^*\nu_{ f}=\eigenv_{ f}\nu_{ f}$ for some $\eigenv_f>0$. Moreover, $\eigenv_f$ is uniquely determined by this equation and satisfies $\eigenv_f=\exp(P(f))$. Here $P\colon\mathcal C( E_A^{\infty})\to\mathbb R$ denotes the \emph{topological pressure function}, which for $\psi\in\mathcal C( E_A^{\infty})$ is defined by
\begin{equation}\label{eq:pressuredef}
P(\psi)\defeq\lim_{n\to\infty}\frac{1}{n}\ln\sum_{\om\in E_A^n}\exp\sup_{\omneu\in[\om]} S_n\psi(\omneu),
\end{equation}
(see \cite[Lem.~1.20]{Bowen_equilibrium}), where we recall that $[\om]\defeq\{\omneu\in E_A^{\infty}\mid \omneu_i=\om_i\ \text{for}\ 1\leq i\leq n(\om)\}$ denotes the $\om$-cylinder set and where the \emph{$n$-th ergodic sum} of a map $f\colon E_A^{\infty}\to\mathbb R$ and $n\in\mathbb N$ is defined to be 
\[
	S_n f\defeq\sum_{k=0}^{n-1} f\circ\sigma^k\quad\text{and}\quad S_0 f\defeq 0.
\]

Further, there exists a unique strictly positive eigenfunction $h_{f}\in\mathcal{C}(E_A^{\infty})$ of $\mathcal L_f$ satisfying $\mathcal L_f h_f=\eigenv_f h_f$ and $\int h_{ f}\textup{d}\nu_{ f}=1$. By $\mu_f$ we denote the $\sigma$-invariant probability measure defined by $\textup{d}\mu_{ f} / \textup{d}\nu_{ f} = h_{ f}$. This is the unique $\sigma$-invariant Gibbs measure for the potential function $f$. 
Additionally, under some normalisation assumptions the iterates of the Perron-Frobenius operator converge to the projection onto the one-dimensional subspace generated by its eigenfunction $\eigenf_f$. To be more precise we have
\begin{equation}\label{eq:convergenceperron}
	\lim_{m\to\infty}\|\eigenv_f^{-m}\mathcal L_f^m\psi - \textstyle{\int}\psi\textup{d}\nu_f\cdot \eigenf_f\|=0\quad\text{for all}\ \psi\in\mathcal C( E_A^{\infty}),
\end{equation}
where $\|\cdot\|$ denotes the supremum norm on $\mathcal C( E_A^{\infty})$.
	The results on the Perron-Frobenius operator quoted above originate mainly from the work of Ruelle, see e.\,g.\  \cite{Ruelle_gas}.

For the geometric potential function $\xi\in\mathcal{C}( E_A^{\infty})$ it can be shown that the \emph{measure theoretical entropy} $\entro(\mu_{-\mdim\xi})$ of the shift-map $\sigma$ with respect to $\mu_{-\mdim \xi}$ is given by
\begin{equation}\label{eq:entropy}
	\entro(\mu_{-\mdim\xi})=\mdim\int\xi\textup{d}\mu_{-\mdim \xi},
\end{equation}
where $\mdim$ denotes the Minkowski dimension of $F$.
This observation follows e.\,g.\ from the variational principle (see \cite[Thm.~1.22]{Bowen_equilibrium}) and the following result, which follows by combining \cite[Thms.~4.2.9, 4.2.11 and 4.2.13]{Urbanski_Buch}.

\begin{proposition}\label{thBedford}
	The Minkowski as well as the Hausdorff dimension of $F$ is equal to the unique real number $t>0$ for which $P(-t\xi)=0$, where $P$ denotes the topological pressure function. 
\end{proposition}

\subsection{Properties of cGDS}\label{sec:selfconfshift}

\begin{proposition}\label{prop:Lvnonempty}
	Let $F$ denote the limit set of a cGDS with aperiodic and irreducible incidence matrix (see Def.~\ref{defn:incidence}). If $F$ satisfies $\leb^1(F)=0$, then $\bigcup_{v\in V}L^v\neq\varnothing$, where $L^v$ is defined in \eqref{eq:Lv}.
\end{proposition}
\begin{proof}
	Assume that $\leb^1(F)=0$ and $\bigcup_{v\in V}L^v=\varnothing$. Then $\bigcup_{v\in V}\left\langle\bigcup_{e\in I_v}\bij[e]\right\rangle=\bigcup_{v\in V}\bigcup_{e\in I_v}\left\langle\bij[e]\right\rangle$. This implies 
	\begin{align*}
		\Phi\left(\bigcup_{v\in V}\left\langle\bigcup_{e\in I_v}\bij[e]\right\rangle\right)
		&= \bigcup_{v\in V}\bigcup_{e'\in T_v}\phi_{e'}\left\langle\bigcup_{e\in I_v}\bij[e]\right\rangle
		=\bigcup_{v\in V}\bigcup_{e'\in T_v}\left\langle\underbrace{\bigcup_{e\in I_v}\bij[e'e]}_{=\bij[e']}\right\rangle\\
		&= \bigcup_{v\in V}\bigcup_{e\in I_v}\left\langle\bij[e]\right\rangle
		=\bigcup_{v\in V}\left\langle\bigcup_{e\in I_v}\bij[e]\right\rangle,
	\end{align*}
	where the second to last equality is due to the fact that the incidence matrix is aperiodic and irreducible. 
	Thus, the set $\bigcup_{v\in V}\left\langle\bigcup_{e\in I_v}\bij[e]\right\rangle$ is invariant under $\Phi$ and hence $F=\bigcup_{v\in V}\left\langle\bigcup_{e\in I_v}\bij[e]\right\rangle$. Since we assume that $\leb^1(F)=0$ and since the sets $\left\langle\bigcup_{e\in I_v}\bij[e]\right\rangle$ are compact non-empty intervals, it follows that $\left\langle\bigcup_{e\in I_v}\bij[e]\right\rangle$ is a singleton for every $v\in V$. 
	Therefore, the cardinality of $F$ is finite which contradicts the fact that the Minkowski dimension of $F$ is positive (see Prop.~\ref{thBedford}).
 \end{proof}

One key property of a cGDS is the bounded distortion property. The following bounded distortion lemma has been obtained in \cite[Lem.~3.2]{KesKom} in the setting of cIFS. The proof follows along the same lines for cGDS giving the following lemma. 

\begin{lemma}[Bounded Distortion]\label{lem:bd}
  There exists a sequence $(\bd_n)_{n\in\mathbb N}$ with $\bd_n>0$ for all $n\in\mathbb N$ and $\lim_{n\to\infty}\bd_n=1$ such that for all $\om,\omneu\in E_A^*$ with $\omneu\om\in E_A^*$ and $x,y\in\phi_{\om}(X_{t(\om_{n(\om)})})$ we have that
  \begin{equation*}
    \bd_{n(\om)}^{-1}\leq\frac{\lvert\phi_{\omneu}'(x)\rvert}{\lvert\phi_{\omneu}'(y)\rvert}\leq\bd_{n(\om)}.
  \end{equation*}
\end{lemma}

\section{Proofs of Thm.~\ref{thm:curvatureresult} and Prop.~\ref{thm:conformalMinkowski}}\label{sec:proofs}

Thm.~\ref{thm:curvatureresult} and Prop.~\ref{thm:conformalMinkowski} are generalisations of \cite[Thms.~2.11 and 2.12]{KesKom} respectively. In the following, we recall the important steps from \cite{KesKom} and point out the necessary modifications. 
The key idea is to prove the statements of Thm.~\ref{thm:curvatureresult} and Prop.~\ref{thm:conformalMinkowski} for the local fractal Euler characteristic and then to apply statements from \cite{RatajWinter} to deduce the respective results concerning the local Minkowski content. 

Without loss of generality we assume that $\{0,1\}\subset F\subseteq[0,1]$ as otherwise the result follows by rescaling. Fix an $\eps>0$ and consider the expression $\leb^0(\partial F_{\eps}\cap(-\infty,\bb])/2$ for some $\bbb\in\mathbb R$.  As in \cite{KesKom} we express $\leb^0(\partial F_{\eps}\cap(-\infty,\bb])/2$ in terms of the image gaps but obtain a different representation because of the non-allowed transitions. 
\begin{equation}\label{eq:lebesgue}
  \frac{\leb^0\big(\partial F_{\eps}\cap(-\infty,\bb]\big)}{2}
    =\underbrace{\sum_{v\in V}\sum_{j=1}^{n_v}\card\{\om\in T_v^*\mid L_{\om}^{v,j}\subseteq(-\infty,\bb],\ \lvert L_{\om}^{v,j}\rvert> 2\eps\}}_{\mbox{\normalsize $\eqdef \Xi(\eps)$}}+\frac{c_1}{2},
\end{equation}
where $c_1\in\{1,2,3\}$ depends on the value of $\bb$. For finding appropriate bounds on $\Xi(\eps)$, we choose an $m\in\mathbb N\cup\{0\}$ such that all \main\ gaps $\{L_{\om}^{v,j}\mid v\in V,j\in\{1,\ldots,n_v\}, \om\in T_v^m\}$ of level $m$ are greater than $2\eps$. For $v\in V$, $j\in\{1,\ldots,n_v\}$ and $\om\in T_v^m$ define
\begin{equation*}
	\Xi_{\om}^{v,j}(\eps)\defeq\card\{\omneu\in T_{i(\om_1)}^*\mid L_{\omneu\om}^{v,j}\subseteq(-\infty,\bb],\ \lvert L_{\omneu\om}^{v,j}\rvert> 2\eps\}.
\end{equation*}
We have the following connection.
\begin{equation}\label{eq:Xisum}
	\sum_{v\in V}\sum_{j=1}^{n_v}\sum_{\om\in T_v^m}\Xi_{\om}^{v,j}(\eps)
	\leq \Xi(\eps)
	\leq \sum_{v\in V}\sum_{j=1}^{n_v}\sum_{\om\in T_v^m}\Xi_{\om}^{v,j}(\eps)+ \sum_{v\in V}n_v\sum_{j=0}^{m-1}(\card E)^j.
\end{equation}
Fix $\bbb\in\mathbb R\setminus F$. Then $F\cap (-\infty,\bb]$ can be expressed as a finite union of sets of the form $\bij[\kappa]$, where $\kappa\in E_A^*$. More precisely, there exists a minimal $l\in\mathbb N$ and $\kappa^{(1)},\ldots,\kappa^{(l)}\in E_A^*$ satisfying
\begin{enumerate}
\item $F\cap(-\infty,\bb]=\bigcup_{j=1}^l\bij[\kappa^{(j)}]$ and
\item $\bij[\kappa^{(i)}]\cap \bij[\kappa^{(j)}]$ contains at most one point for all $i\neq j\in\{1,\ldots,l\}$.
\end{enumerate}
Then for $\kappa\defeq \bigcup_{j=1}^l[\kappa^{(j)}]$ the function $\mathds 1_{\kappa}$ is H\"older continuous. Fix $\om\in E_A^m$.
Using the bounded distortion constant $\bd_{n(\om)}$ of $\Phi$ on $\phi_{\om}(X_{t(\om_{n(\om)})})$ (see Lem.~\ref{lem:bd}), we can provide upper bounds for $\Xi_{\om}^{v,j}(\eps)$, namely for an arbitrary $\om^v\in  I_v^{\infty}$ we have

\begin{align}
	\Xi_{\om}^{v,j}(\eps)
	&\leq \sum_{n=0}^{\infty}\sum_{\omneu\in T_{i(\om_1)}^n}\mathds 1_{\kappa}(\omneu\om \om^v)\mathds 1_{\big{\{}\lvert\phi_{\omneu}'(\bij\om \om^v)\rvert\cdot\bd_{n(\om)}\cdot\lvert L_{\om}^{v,j}\rvert>2\eps\big{\}}}+ \overline{c}_2(\om^v,\kappa)\nonumber\\
	&\leq \underbrace{\sum_{n=0}^{\infty}\sum_{\omneu\in T_{i(\om_1)}^n}\mathds 1_{\kappa}(\omneu\om \om^v)\mathds 1_{\big{\{}\lvert\phi_{\omneu}'(\bij\om \om^v)\rvert\cdot\bd_{n(\om)}\cdot\lvert L_{\om}^{v,j}\rvert\geq 2\eps\big{\}}}}_{\mbox{\normalsize $\eqdef\overline{A}_{\om}^{v,j}(\om^v,\eps,\kappa)$}}+ \overline{c}_2(\om^v,\kappa),\label{eq:Aomover}
\end{align}
where the constant $\overline{c}_2(\om^v,\kappa)$ is needed, since $L_{\omneu\om}^{v,j}\subseteq(-\infty,\bb]$ does not necessarily imply $\omneu\om \om^v\in \kappa$ for an arbitrary $\om^v\in I_v^{\infty}$. However, if $n(\omneu)\geq\max_{j=1,\ldots,l}n(\kappa^{(j)})$, either $[\omneu\om]\subseteq \kappa$ or $[\omneu\om]\cap \kappa=\varnothing$. Hence, there are only finitely many $\omneu\in T_{i(\om_1)}^*$ for which $L_{\omneu\om}^{v,j}\subseteq(-\infty,\bb]$ does not imply $\omneu\om \om^v\in \kappa$ for all $\om^v\in I_v^{\infty}$. Letting $\overline{c}_2(\om^v,\kappa)\in\mathbb R$ denote this finite number shows \eqref{eq:Aomover} for all $\eps>0$. Likewise, there exists a constant $\underline{c}_2(\om^v,\kappa)\in\mathbb R$ such that for all $\eps>0$
\begin{equation*}
	\Xi_{\om}^{v,j}(\eps)
	\geq \sum_{n=0}^{\infty}\sum_{\omneu\in T_{i(\om_1)}^n}\mathds 1_{\kappa}(\omneu\om \om^v)\mathds 1_{\big{\{}\lvert\phi_{\omneu}'(\bij\om \om^v)\rvert\cdot\bd_{n(\om)}^{-1}\cdot\lvert L_{\om}^{v,j}\rvert>2\eps\big{\}}}- \underline{c}_2(\om^v,\kappa).
\end{equation*}
It follows that for all $\beta>1$ we have that 
\begin{equation}\label{eq:Aomunder}
	\Xi_{\om}^{v,j}(\eps)
	\geq \underbrace{\sum_{n=0}^{\infty}\sum_{\omneu\in T_{i(\om_1)}^n}\mathds 1_{\kappa}(\omneu\om \om^v)\mathds 1_{\big{\{}\lvert\phi_{\omneu}'(\bij\om \om^v)\rvert\cdot\bd_{n(\om)}^{-1}\cdot\lvert L_{\om}^{v,j}\rvert\geq 2\eps\beta\big{\}}}}_{\mbox{\normalsize $\eqdef\underline{A}_{\om}^{v,j}(\om^v,\eps\beta,\kappa)$}}- \underline{c}_2(\om^v,\kappa).
\end{equation}
For every $v\in V$ fix an $\om^v\in I_v^{\infty}$.
Combining \eqref{eq:lebesgue} to \eqref{eq:Aomunder} implies for all $m\in\mathbb N$ and all $\beta>1$ that 
\begin{align}
  \overline{C}_0^f(F,(-\infty,\bb])
    &\leq \limsup_{\eps\to 0}\eps^{\mdim}\sum_{v\in V}\sum_{j=1}^{n_v}\sum_{\om\in T_v^m} \overline{A}_{\om}^{v,j}(\om^v,\eps,\kappa)\quad\text{and}\label{eq:fcmoben}\\
    \underline{C}_0^f(F,(-\infty,\bb])
      &\geq \liminf_{\eps\to 0}\eps^{\mdim}\sum_{v\in V}\sum_{j=1}^{n_v}\sum_{\om\in T_v^m} \underline{A}_{\om}^{v,j}(\om^v,\eps\beta,\kappa).\label{eq:fcmunten}	
\end{align}
The next step in the proofs is to apply renewal theorems by Lalley \cite{Lalley} and slight extensions by the authors \cite{KesKom}, in order to obtain asymptotics for both expressions $\overline{A}_{\om}^{v,j}(\om^v,\eps,\kappa)$ and $\underline{A}_{\om}^{v,j}(\om^v,\eps\beta,\kappa)$.
For applying the renewal theorems, note that
\begin{align}
  &\sum_{\omneu\in T_{i(\om_1)}^n}\mathds 1_{\kappa}(\omneu\om \om^v)\cdot\mathds 1_{\big{\{}\lvert\phi_{\omneu}'(\bij\om \om^v)\rvert\cdot\bd_{n(\om)}^{\pm 1}\cdot\lvert L_{\om}^{v,j}\rvert\geq2\eps\big{\}}}\nonumber\\
  &\qquad = \sum_{\omneu:\sigma^n \omneu=\om \om^v}\mathds 1_{\kappa}(\omneu)\cdot\mathds 1_{\Bigg{\{}
{\displaystyle\sum_{k=1}^{n}}-\ln\lvert\phi_{\omneu_k}'(\bij\sigma^k \omneu)\rvert\leq-\ln\frac{2\eps}{\lvert L_{\om}^{v,j}\rvert\bd_{n(\om)}^{\pm 1}}\Bigg{\}}}
\nonumber\\
  &\qquad= \sum_{\omneu:\sigma^n \omneu=\om \om^v}\mathds 1_{\kappa}(\omneu)\cdot\mathds 1_{\bigg{\{}S_n\xi(\omneu)\leq-\ln\frac{2\eps}{\lvert L_{\om}^{v,j}\rvert\bd_{n(\om)}^{\pm 1}}\bigg{\}}}.\label{eq:Snpm}
\end{align}
Moreover, note that the hypotheses and Rem.~\ref{rmk:alpha} imply that the geometric potential function $\xi$ is H\"older continuous and strictly positive. The unique $s>0$ for which $\eigenv_{-s\xi}=1$ is precisely the Minkowski dimension $\mdim$ of $F$, which results by combining the fact that $\eigenv_{-s\xi}=\exp(P(-s\xi))$ for each $s>0$ and Prop.~\ref{thBedford}.
We will insert the asymptotics for $\overline{A}_{\om}^{v,j}(\om^v,\eps,\kappa)$ and $\underline{A}_{\om}^{v,j}(\om^v,\eps\beta,\kappa)$, that the renewal theorems yield, into \eqref{eq:fcmoben} and \eqref{eq:fcmunten}. In this way we will obtain an upper bound for $ \overline{C}_0^f(F,(-\infty,\bb])$ and a lower bound for $ \underline{C}_0^f(F,(-\infty,\bb])$. For deducing statements on $ \overline{C}_0^f(F,(-\infty,\bb])$ and $ \underline{C}_0^f(F,(-\infty,\bb])$ from these bounds we need the following lemma, which is an adaptation of \cite[Lem.~4.1]{KesKom}.

\begin{lemma}\label{lem:Upsilon}
  For every $v\in V$ fix an $\om^v\in I_v^{\infty}$. Then for an arbitrary $\Upsilon\in\mathbb R$ we have that
  \begin{enumerate}
  \item
    $\Upsilon\leq\sum_{v\in V}\sum_{j=1}^{n_v}\sum_{\om\in T_v^m}\eigenf_{-\mdim\xi}(\om \om^v)\big(\lvert L_{\om}^{v,j}\rvert\bd_{m}\big)^{\mdim}$ for all $m\in\mathbb N$ implies
    \[
    \Upsilon\leq\liminf_{m\to\infty}\sum_{v\in V}\sum_{j=1}^{n_v}\sum_{\om\in T_v^m}\lvert L_{\om}^{v,j}\rvert^{\mdim}.
    \]
  \item 
    $\Upsilon\geq\sum_{v\in V}\sum_{j=1}^{n_v}\sum_{\om\in T_v^m}\eigenf_{-\mdim\xi}(\om \om^v)\big(\lvert L_{\om}^{v,j}\rvert\bd_{m}^{-1}\big)^{\mdim}$ for all $m\in\mathbb N$ implies
    \[
    \Upsilon\geq\limsup_{m\to\infty}\sum_{v\in V}\sum_{j=1}^{n_v}\sum_{\om\in T_v^m}\lvert L_{\om}^{v,j}\rvert^{\mdim}.
    \]
  \end{enumerate}
\end{lemma}

\begin{proof}
First, we approximate the eigenfunction $\eigenf_{-\mdim\xi}$ of the Perron-Frobenius operator $\mathcal L_{-\mdim\xi}$. By induction it follows that
\[
\mathcal L_{-\mdim\xi}^n \mathbf{1}(\om)=\sum_{\omneu\in T_{i(\om_1)}^n}\lvert \phi_{\omneu}'(\bij \om)\rvert^{\mdim}
\]
for each $\om\in E_A^{\infty}$ and $n\in\mathbb N$, where $\mathbf{1}$ is the constant one-function on $E_A^{\infty}$. 
$\mathcal L_{-\mdim\xi}^n \mathbf{1}$ converges uniformly to the eigenfunction $\eigenf_{-\mdim\xi}$ when taking $n\to\infty$ by \eqref{eq:convergenceperron}. Thus, 
\[
\forall t>0\ \exists M\in\mathbb N\colon\forall n\geq M,\ \forall\, \om\in E_A^{\infty}\colon \bigg{\lvert}\sum_{\omneu\in T_{i(\om_1)}^n}\lvert\phi_{\omneu}'(\bij \om)\rvert^{\mdim}-\eigenf_{-\mdim\xi}(\om)\bigg{\rvert}<t.
\]
Furthermore, the important bounded distortion lemma (Lem.~\ref{lem:bd}) states that 
\[
\forall t'>0\ \exists M'\in\mathbb N\colon\forall m\geq M'\colon \lvert\bd_{m}-1\rvert<t'.
\]
Thus, for all $n\geq M$ and $m\geq M'$
\begin{align*}
  \Upsilon
  &\leq\sum_{v\in V}\sum_{j=1}^{n_v}\sum_{\om\in T_v^m}\eigenf_{-\mdim\xi}(\om \om^v)\left(\lvert L_{\om}^{v,j}\rvert\bd_{m}\right)^{\mdim}\\
  &\leq \sum_{v\in V}\sum_{j=1}^{n_v}\sum_{\om\in T_v^m}\left(\sum_{\omneu\in T_{i(\om_1)}^n}\lvert\phi_{\omneu}'(\bij\om \om^v)\rvert^{\mdim}+t\right)\lvert L_{\om}^{v,j}\rvert^{\mdim}(1+t')^{\mdim}\\
  &\leq \sum_{v\in V}\sum_{j=1}^{n_v}\sum_{\om\in T_v^m}\sum_{\omneu\in T_{i(\om_1)}^n}\lvert L_{\omneu\om}^{v,j}\rvert^{\mdim}(1+t')^{2\mdim}+ t(1+t')^{\mdim}\sum_{v\in V}\sum_{j=1}^{n_v}\sum_{\om\in T_v^m}\lvert L_{\om}^{v,j}\rvert^{\mdim}\\
  &\eqdef A_{m,n}.
\end{align*}
Hence, for all $t,t'>0$,
\begin{align*}
	\Upsilon
	&\leq \liminf_{m\to\infty}\liminf_{n\to\infty} A_{m,n}\\
	&\leq \left(1+t'\right)^{2\mdim}\liminf_{m\to\infty}\liminf_{n\to\infty}\sum_{v\in V}\sum_{j=1}^{n_v}\sum_{\om\in T_v^{m}}\sum_{\omneu\in T_{i(\om_1)}^n}\lvert L_{\omneu\om}^{v,j}\rvert^{\mdim}\\
	&\hphantom{\leq} +t(1+t')^{\mdim}\limsup_{m\to\infty}\sum_{v\in V}\sum_{j=1}^{n_v}\sum_{\om\in T_v^m}\lvert L_{\om}^{v,j}\rvert^{\mdim}.
\end{align*}
We have that $\sum_{v\in V}\sum_{j=1}^{n_v}\sum_{\om\in T_v^m}\lvert L_{\om}^{v,j}\rvert^{\mdim}\leq\sum_{v\in V}n_v\cdot\sum_{\om\in T_v^m}\|\phi'_{\om}\|^{\mdim}\eqdef a_m$, where $\|\cdot\|$ denotes the supremum norm. The assertion follows by letting $t$ and $t'$ tend to zero, since the sequence $(a_m)_{m\in\mathbb N}$ is bounded by \cite[Lem.~4.2.12]{Urbanski_Buch} together with Rem.~\ref{rmk:finitelyprimitive}.
Analogously, the lower bound in the second case can be proved.
 \end{proof}

\subsection{The Non-lattice Case}
\begin{proof}[Proof of Thm.~\ref{thm:curvatureresult}\ref{curvatureresult:nonlattice}]
Let us fix the notation from the beginning of Sec.~\ref{sec:proofs}. 

If $\mathds 1_{\kappa}$ is identically zero, then $C_0^f(F,(-\infty,\bb])=0=\nu(F\cap(-\infty,\bb])$. If $\mathds 1_{\kappa}$ is not identically zero then combining \eqref{eq:Aomover}, \eqref{eq:Aomunder} and \eqref{eq:Snpm} with the fact that $\mathds 1_{\kappa}$ is H\"older continuous allows us to apply Lalley's renewal theorem \cite[Thm.~1]{Lalley} (see also \cite[Prop.~3.8]{KesKom}, where the theorem is stated using our notation, but for the case that $E_A^{\infty}=E^{\mathbb N}\eqdef \Sigma^{\infty}$) to  $\overline{A}_{\om}^{v,j}(\om^v,\eps,\kappa)$ and $\underline{A}_{\om}^{v,j}(\om^v,\eps\beta,\kappa)$, where $v\in V$, $j\in\{1,\ldots,n_v\}$, $\om\in T_v^*$, $\om^v\in I_v^{\infty}$ and $\beta> 1$. This gives the following asymptotics.

\begin{align}
  \overline{A}_{\om}^{v,j}(\om^v,\eps,\kappa)
  &\sim \frac{\nu_{-\mdim\xi}(\kappa)}{\mdim\int\xi\textup{d}\mu_{-\mdim\xi}}\cdot\eigenf_{-\mdim\xi}(\om \om^v)\cdot(2\eps)^{-\mdim}\big(\lvert L_{\om}^{v,j}\rvert\bd_{n(\om)}\big)^{\mdim},\label{eq:Aomoverasym}\\
  \underline{A}_{\om}^{v,j}(\om^v,\eps\beta,\kappa)
  &\sim  \frac{\nu_{-\mdim\xi}(\kappa)}{\mdim\int\xi\textup{d}\mu_{-\mdim\xi}}
  \cdot\eigenf_{-\mdim\xi}(\om \om^v)\cdot(2\eps\beta)^{-\mdim}\big(\lvert L_{\om}^{v,j}\rvert\bd_{n(\om)}^{-1}\big)^{\mdim}\label{eq:Aomunderasym}
\end{align}
as $\eps\to 0$ uniformly for $\om^v\in I_v^{\infty}$.  On combining \eqref{eq:fcmoben}, \eqref{eq:Aomoverasym} (resp.\   \eqref{eq:fcmunten}, \eqref{eq:Aomunderasym}) with Lem.~\ref{lem:Upsilon} we obtain in a similar way to \cite[proof of Thm.~2.11(ii)]{KesKom} that $\overline{C}_{0}^{f}(F,(-\infty,b])=\underline{C}_{0}^{f}(F,(-\infty,b])$ and thus that
\[
	C_0^f(F,(-\infty,\bb])
=\frac{2^{-\mdim}}{\entro(\mu_{-\mdim\xi})}\lim_{m\to\infty}\sum_{v\in V}\sum_{j=1}^{n_v}\sum_{\om\in T_v^m}\lvert L_{\om}^{v,j}\rvert^{\mdim}\cdot \nu(F\cap(-\infty,\bb])
\]
holds for every $\bbb\in\mathbb R\setminus F$. As $\mathbb R\setminus F$ is dense in $\mathbb R$ the assertion concerning the local fractal Euler characteristic follows.
The result on the local Minkowski content now follows by applying \cite[Cor.~3.2]{RatajWinter} (see also \cite[Thm.~3.13]{KesKom}), as for every $\bbb\in\mathbb R\setminus F$ we have that $F_{\eps}\cap(-\infty,\bb]=\left(F\cap(-\infty,\bb]\right)_{\eps}$ for sufficiently small $\eps>0$.
 \end{proof}

\subsection{The Lattice Case}
This subsection addresses Thm.~\ref{thm:curvatureresult}\ref{curvatureresult:lattice} and Prop.~\ref{thm:conformalMinkowski}.

\begin{proof}[Proof of Thm.~\ref{thm:curvatureresult}\ref{curvatureresult:lattice}]
The statement that neither the local Minkowski content nor the local fractal Euler characteristic exists if the maps $\phi_e$ are all analytic, is a direct consequence of Thm.~\ref{thm:imafcm}\ref{imafcm_lattice} together with Thm.~\ref{thm:analytic}, which will both be proven in Sec.~\ref{sec:specialcases}. Therefore, we now focus on the boundedness and positivity. For this, we proceed as in the proof of \cite[Thm.~2.11(iii)]{KesKom}. We fix the notation from the beginning of Sec.~\ref{sec:proofs}.
As $\xi$ is lattice, there are $\ze,\psi\in\mathcal C( E_A^{\infty})$ such that
\[
	\xi-\ze=\psi-\psi\circ \sigma,
\]
where the range of $\ze$ is contained in a discrete subgroup of $\mathbb R$. Let $\aaa>0$ be the maximal real number for which $\ze( E_A^{\infty})\subseteq\aaa\mathbb Z$. Recall from the beginning of Sec.~\ref{sec:proofs} that the hypotheses and Rem.~\ref{rmk:alpha} imply that $\xi$ is H\"older continuous and strictly positive and that the unique $s>0$ for which $\eigenv_{-s\xi}=1$ is the Minkowski dimension $\mdim$ of $F$ (see Prop.~\ref{thBedford}). Note that $\mathds 1_{\kappa}$ is H\"older continuous and that we can assume that $\mathds 1_{\kappa}$ is not identically zero. Combining \eqref{eq:Aomover}, \eqref{eq:Aomunder} and \eqref{eq:Snpm}, we see that we can apply the extended version of Lalley's renewal theorem \cite[Thm.~3]{Lalley} given in \cite[Thm.~3.9]{KesKom} to $\overline{A}_{\om}^{v,j}(\om^v,\eps,\kappa)$ and $\underline{A}_{\om}^{v,j}(\om^v,\eps\beta,\kappa)$. This yields the following asymptotics.
\begin{align}
  \overline{A}_{\om}^{v,j}(\om^v,\eps,\kappa)
  &\sim U_{\om}(\om^v)\int_{\kappa} \ee^{-\mdim\aaa\left\lceil\frac{\psi(y)-\psi(\om \om^v)}{\aaa}+\frac{1}{\aaa}\ln\frac{2\eps}{\lvert L_{\om}^{v,j}\rvert\bd_{n(\om)}}\right\rceil}\textup{d}\nu_{-\mdim\ze}(y),\label{eq:Aomoverasym:lattice}\\
  \underline{A}_{\om}^{v,j}(\om^v,\eps\beta,\kappa)
  &\sim U_{\om}(\om^v)\int_{\kappa} \ee^{-\mdim\aaa\left\lceil\frac{\psi(y)-\psi(\om \om^v)}{\aaa}+\frac{1}{\aaa}\ln\frac{2\eps\beta\bd_{n(\om)}}{\lvert L_{\om}^{v,j}\rvert}\right\rceil}\textup{d}\nu_{-\mdim\ze}(y)\label{eq:Aomunderasym:lattice}
\end{align}
as $\eps\to 0$ uniformly for $\om^v\in I_v^{\infty}$, where
\begin{equation}\label{eq:W}
  U_{\om}(\om^v)\defeq \frac{\aaa\eigenf_{-\mdim\ze}(\om \om^v)}{(1-\ee^{-\mdim\aaa})\int\ze\textup{d}\mu_{-\mdim\ze}}.
\end{equation}
Consecutively applying \eqref{eq:fcmoben}, \eqref{eq:Aomoverasym:lattice} and using that $\ee^{-\lfloor x\rfloor}\leq \ee^{-x}$ for any $x\in\mathbb R$ we obtain with Lem.~\ref{lem:Upsilon} (as in \cite[p.\,2498]{KesKom}) that
  \begin{equation*}
    \overline{C}_0^f(F,\mathbb R)
    \leq \liminf_{m\to\infty}\sum_{v\in V}\sum_{j=1}^{n_v}\sum_{\om\in T_v^m}
		\lvert L_{\om}^{v,j}\rvert^{\mdim}\frac{\aaa 2^{-\mdim}}{(1-\ee^{-\mdim\aaa})\int\ze\textup{d}\mu_{-\mdim\ze}}\eqdef c_0.
  \end{equation*}
  The number $c_0$ is positive and finite because 
\[\sum_{v\in V}\sum_{j=1}^{n_v}\sum_{\om\in T_v^m}\lvert L_{\om}^{v,j}\rvert^{\mdim}\leq \sum_{v\in V}\sum_{j=1}^{n_v}\sum_{\om\in T_v^m}\|\phi'_{\om}\|^{\mdim}\eqdef a_m,\]
 where $\|\cdot\|$ denotes the supremum norm, and the sequence $(a_m)_{m\in\mathbb N}$ is bounded by \cite[Lem.~4.2.12]{Urbanski_Buch}.
  
  That $\underline{C}_0^f(F,\mathbb R)$ is positive can be concluded from \eqref{eq:fcmunten}, \eqref{eq:Aomunderasym:lattice} and Lem.~\ref{lem:Upsilon} in the same way (see also \cite[p.\,2499]{KesKom}).

  The results on $\underline{C}_1^f(F,B)$ and $\overline{C}_1^f(F,B)$ are now straightforward applications of \cite[Cor.~3.2]{RatajWinter} (see also \cite[Thm.~3.13]{KesKom}).
 \end{proof}

\begin{proof}[Proof of Prop.~\ref{thm:conformalMinkowski}]
	Equation \eqref{eq:Mclattice} follows from Thm.~\ref{thm:curvatureresult}\ref{curvatureresult:lattice}. For the second statement, we use \cite[Lem.~3.12]{KesKom}, which in \cite{KesKom} was proven for the case that $E_A^{\infty}=E^{\mathbb N}\eqdef \Sigma^{\infty}$, but the same proof works in the present more general situation. The hypotheses and \cite[Lem.~3.12]{KesKom} together imply that for every $v\in V$, $j\in\{1,\ldots, n_v\}$, $\om^v\in I_v^{\infty}$ and $\om\in T_v^m$ we have that
  \begin{align*}
    \overline{U}
    &\defeq\lim_{\eps\to 0}\eps^{\mdim}\int\ee^{-\mdim\aaa\left\lceil\frac{\psi(y)-\psi(\om \om^v)}{\aaa}+\frac{1}{\aaa}\ln\frac{2\eps}{\lvert L_{\om}^{v,j}\rvert\bd_m}\right\rceil}\textup{d}\nu_{-\mdim\ze}(y)\cdot \left(\frac{2}{\lvert L_{\om}^{v,j}\rvert\bd_m}\right)^{\mdim}\quad\text{and}\\
    \underline{U}
    &\defeq\lim_{\eps\to 0}\eps^{\mdim}\int\ee^{-\mdim\aaa\left\lceil\frac{\psi(y)-\psi(\om \om^v)}{\aaa}+\frac{1}{\aaa}\ln\frac{2\eps\bd_m}{\lvert L_{\om}^{v,j}\rvert}\right\rceil}\textup{d}\nu_{-\mdim\ze}(y)\cdot \left(\frac{2\bd_m}{\lvert L_{\om}^{v,j}\rvert}\right)^{\mdim}
  \end{align*}
	are independent of $\om$, $v$ and $j$ and are equal, i.\,e.\ $\overline{U}=\underline{U}\eqdef U$.
	Combining \eqref{eq:fcmoben} with \eqref{eq:Aomoverasym:lattice} and \eqref{eq:fcmunten} with \eqref{eq:Aomunderasym:lattice} and applying Lem.~\ref{lem:Upsilon} gives $\overline{C}_0^f(F,\mathbb R)=\underline{C}_0^f(F,\mathbb R)$. These steps are carried out in the case $E_A^{\infty}=E^{\mathbb N}$ in \cite[p.\,2499 f.]{KesKom} in detail. An application of \cite[Cor.~3.2]{RatajWinter} (see also \cite[Thm.~3.13]{KesKom}) then completes the proof.
 \end{proof}

\subsection{Average Quantities}
\begin{proof}[Proof of Thm.~\ref{thm:curvatureresult}\ref{curvatureresult:average}]
	This statement follows in direct analogy to the proof of \cite[Thm.~2.11(i)]{KesKom}. Slight modifications to match the present setting are the following. For $m\in\mathbb N$ set ${M\defeq \min\{\lvert L_{\om}^{v,j}\rvert\mid v\in V, j\in\{1,\ldots,n_v\}, \om\in T_v^m\}/2}$. 
Note that the analogues of \cite[(4.1)\,--\,(4.3)]{KesKom} are \eqref{eq:lebesgue}\,--\,\eqref{eq:Aomover} and that $\sum_{i=1}^L\sum_{\om\in\Sigma^m}\overline{A}_{\om}^i(x,\eps,\kappa)$ needs to be replaced by 
\[
	\sum_{v\in V}\sum_{j=1}^{n_v}\sum_{\om\in T_v^m} \overline{A}_{\om}^{v,j}(\om^v,\eps,\kappa).
\]
The analogue of \cite[(4.7)]{KesKom} is \eqref{eq:Snpm}, the analogue of \cite[Lem.~4.1]{KesKom} is Lem.~\ref{lem:Upsilon} and the analogue of \cite[Thm.~2.11(iii)]{KesKom} is Thm.~\ref{thm:curvatureresult}\ref{curvatureresult:lattice}.
\end{proof}

\section{Proofs of Thms.~\ref{thm:similars} to \ref{thm:analytic}}\label{sec:specialcases}

Here, we provide the proofs of the results concerning limit sets of sGDS (Thm.~\ref{thm:similars}) and piecewise $\mathcal C^{1+\alpha}$-diffeomorphic images of limit sets of sGDS (Thms.~\ref{thm:imafcm}, \ref{thm:analytic}). 

\subsection{sGDS}

\renewcommand{\theenumi}{(\Alph{enumi})}

\begin{proof}[Proof of Thm.~\ref{thm:similars}]
  Throughout the proof let $\Phi\defeq(\phi_e)_{e\in E}$ denote an sGDS, meaning that $\phi_e$ is a similarity for every $e\in E$. Let $r_e\in(0,1)$ denote a Lipschitz constant of $\phi_e$ for $e\in E$. Further, set $r_{\om}\defeq r_{\om_1}\cdots r_{\om_n}$ for a finite word $\om=\om_1\cdots\om_n\in E_A^n$.
  Part \ref{ss:average} follows from Thm.~\ref{thm:curvatureresult}\ref{curvatureresult:average} with the following considerations. For $v\in V$ and $\om^v\in I_v^{\infty}$, \eqref{eq:convergenceperron} implies 
  \begin{equation}\label{eq:eigenfapproxi}
    \lim_{m\to\infty}\sum_{\om\in T_v^m} r_{\om}^{\mdim}
    =\lim_{m\to\infty}\mathcal L_{-\mdim\xi}^m \mathbf{1}(\om^v)
    =\eigenf_{-\mdim\xi}(\om^v)
  \end{equation}
with $\mathbf 1$ denoting the constant one-function on $E_A^{\infty}$.
  Moreover, since $\phi_e$ are similarities, $\lvert L_{\om}^{v,j}\rvert=r_{\om}\lvert L^{v,j}\rvert$ for all $v\in V$, $\om\in T_v^*$ and $j\in\{1,\ldots, n_v\}$.
  Thus, $c$ from \eqref{eq:constantthm} simplifies to 
  \begin{equation}
    c=\sum_{v\in V}\sum_{j=1}^{n_v}\eigenf_{-\mdim\xi}(\om^v)\lvert L^{v,j}\rvert^{\mdim},
  \end{equation}
   showing the assertion of \ref{ss:average}.
  
  In order to prove \ref{ss:lattice} we proceed as in the proof of \cite[Thm.~2.14]{KesKom}, where the statement is shown for self-similar sets. A crucial discrepancy to \cite{KesKom} is that here $\eigenf_{-\mdim\xi}\equiv 1$ is not necessarily satisfied. 
  With the same arguments as in the proof of \cite[Thm.~2.14]{KesKom} (replacing $E$ and $R_{\om}E$ in \cite{KesKom} by $F$ and $\pi[\om]$ respectively) we see that for fixed $v\in V$ and arbitrary $\om^v\in I_v^{\infty}$ there exists a constant $\widetilde{c}\geq 0$, which depends on the number of sets $\bij[\om]$ whose union is $F\cap B$, such that
  \begin{align}
    \leb^0\big(\partial F_{\ee^{-T}}\cap B \big)/2
    &\stackrel{\eqref{eq:lebesgue}}{=}\sum_{v\in V}\sum_{j=1}^{n_v}\card\{\om\in T_v^*\mid L_{\om}^{v,j}\subseteq B,\ \lvert L_{\om}^{v,j}\rvert> 2\ee^{-T}\}+\widetilde{c}\nonumber\\
    &\ \ =\ \ \sum_{v\in V}\sum_{j=1}^{n_v}\sum_{n=0}^{\infty}\sum_{\om\in T_v^n}\mathds 1_{\kappa}(\om \om^v)\mathds 1_{\big{\{}\lvert\phi'_{\om}(\om^v)\rvert\cdot\lvert L^{v,j}\rvert> 2\ee^{-T}\big{\}}}+\widetilde{c}\nonumber\\
    &\ \ =\ \ \sum_{v\in V}\sum_{j=1}^{n_v}\sum_{n=0}^{\infty}\sum_{\omneu:\sigma^n \omneu=\om^v} \mathds 1_{\kappa} (\omneu)\mathds 1_{\big{\{}S_n\xi(\omneu)<-\ln\frac{2\ee^{-T}}{\lvert L^{v,j}\rvert}\big{\}}}+\widetilde{c}\nonumber\\
    &\ \ \sim\ \ \sum_{v\in V}\sum_{j=1}^{n_v}\frac{\aaa\eigenf_{-\mdim\xi}(\om^v)\nu_{-\mdim\xi}(\kappa)}{(1-\ee^{-\mdim\aaa})\int\xi\textup{d}\mu_{-\mdim\xi}}\cdot \ee^{-\mdim\aaa\left\lceil\aaa^{-1}\ln\frac{2\ee^{-T}}{\lvert L^{v,j}\rvert}\right\rceil}+\widetilde{c}\label{eq:ssasymptoticproof}
  \end{align}
  as $T\to\infty$, where the last asymptotic is obtained by applying \cite[Thm.~3.9 and Rem.~3.10]{KesKom}, a slight extension of \cite[Thm.~3]{Lalley}. As in \cite{KesKom} we introduce the function $f\colon\mathbb R^+\to\mathbb R^+$ which here is given by
  \[
  f(T)\defeq \ee^{-\mdim T}\frac{\aaa\nu(B)}{(1-\ee^{-\mdim\aaa})\entro(\mu_{-\mdim\xi})}\sum_{v\in V}\sum_{j=1}^{n_v} \eigenf_{-\mdim\xi}(\om^v)\ee^{-\mdim\aaa\left\lceil\frac{1}{\aaa}\ln\frac{2\ee^{-T}}{\lvert L^{v,j}\rvert}\right\rceil}.
  \]
  By the asymptotic given in \eqref{eq:ssasymptoticproof}, we know that for all $t>0$ there exists an $M\in\mathbb N$ such that for all $T\geq M$ we have
  \begin{equation*}
    (1-t)\mdim f(T)
    \leq \ee^{-\mdim T}\leb^0(\partial F_{\ee^{-T}}\cap B)/2
    \leq(1+t)\mdim f(T)+c\ee^{-\mdim T}.
  \end{equation*}
  Clearly, $f$ is a strictly positive, bounded and periodic function with period $\aaa$. Moreover, $f$ is piecewise continuous with a finite number of discontinuities in an interval of length $\aaa$. Additionally, on every interval, where $f$ is continuous, $f$ is strictly decreasing. Therefore $f$ is not equal to an almost everywhere constant function. Thus, all conditions of \cite[Lem.~5.1]{KesKom} (whose proof works in exactly the same way for cGDS instead of cIFS) are satisfied which shows the statement.
 \end{proof}

\subsection{Piecewise $\mathcal C^{1+\alpha}$-diffeomorphic Images of Limit Sets of sGDS}
Here, we consider the case that $F$ is the image of the limit set $K$ of an sGDS under a piecewise $\mathcal C^{1+\alpha}$-diffeomorphism. Throughout, we fix the notation from Thm.~\ref{thm:imafcm}.
By definition, each $g_v$ is bi-Lipschitz. Therefore, the Minkowski dimensions of $K$ and $F$ coincide (see e.\,g.\ \cite[Cor.~2.4 and Sec.~3.2]{Falconer_Foundation}) and are both denoted by $\mdim$.

The similarities $(R_e)_{e\in E}$ generating $K$ and the mappings $(\phi_e)_{e\in E}$ generating $F$ are connected through the equations 
\[
	\phi_e=g_{i(e)}\circ R_e\circ g_{t(e)}^{-1}
\]
for each $e\in E$. We denote by $\widetilde{\bij}$ and $\bij$ the natural code maps from $E_A^{\infty}$ to $K$ and $F$ respectively. 
If we further let $(r_e)_{e\in E}$ denote the respective similarity ratios of $(R_e)_{e\in E}$, i.\,e.\ $r_e\defeq\|R'_e\|$, we have the following list of observations.

\renewcommand{\theenumi}{(\Alph{enumi})}
\begin{enumerate}
\item\label{en:diffbar} Each map $\phi_e\colon X_{t(e)}\to X_{i(e)}$ is differentiable with derivative
  \[
  \phi'_e=\frac{g_{i(e)}'\circ R_e\circ g_{t(e)}^{-1}}{g_{t(e)}'\circ g_{t(e)}^{-1}}\cdot r_e.
  \] 
\item\label{en:geopot} The geometric potential function $\ze$ associated with $K$ is given by $\ze(\om)=-\ln r_{\om_1}$;
the geometric potential function $\xi$ associated with $F$ is given by $\xi(\om)=-\ln\lvert g_{t(\om_1)}'(g_{t(\om_1)}^{-1}(\bij\om))\rvert +\ln\lvert g'_{t(\om_2)}(g_{t(\om_2)}^{-1}(\bij \sigma\om))\rvert-\ln r_{\om_1}$, where $\om=\om_1\om_2\cdots\in E_A^{\infty}$.
 Thus $\ze$ is non-lattice, if and only if $\xi$ is non-lattice.
\item\label{en:gibbs} The unique $\sigma$-invariant Gibbs measure for the potential function $-\mdim\xi$ coincides with the unique $\sigma$-invariant Gibbs measure for the potential function $-\mdim\ze$, i.\,e.\ $\mu_{-\mdim\xi}=\mu_{-\mdim\ze}$ (see e.\,g.\ \cite[Thm.~2.2.7]{Urbanski_Buch}). 
\item From \ref{en:geopot} and \ref{en:gibbs} we obtain that
  \[
  \hspace{0.5cm}\entro(\mu_{-\mdim\xi})
  = \int\xi\textup{d}\mu_{-\mdim\xi}
  = \int\ze\textup{d}\mu_{-\mdim\ze}
  =\entro(\mu_{-\mdim\ze}).
  \]
\end{enumerate} 
Further, let $\{\widetilde{L}^{v,j}\}_{v\in V, j\in\{1,\ldots,n_v\}}$ denote the primary gaps of $K$ and denote by $\{\widetilde{L}_{\om}^{v,j}\}_{v\in V, j\in\{1,\ldots,n_v\}}$ the \main\ gaps of $K$ for each $\om\in E_A^*$. Likewise we let $\{L^{v,j}\}_{v\in V, j\in\{1,\ldots,n_v\}}$ and $\{L_{\om}^{v,j}\}_{v\in V, j\in\{1,\ldots,n_v\}}$ respectively denote the primary and the \main\ gaps of $F$.
\begin{enumerate}\setcounter{enumi}{4}
\item\label{en:nu} The $\mdim$-conformal measure $\nu$ associated with $(\phi_e)_{e\in E}$ and the push-forward measure of the $\mdim$-conformal measure $\widetilde{\nu}$ associated with $(R_e)_{e\in E}$ are absolutely continuous with Radon-Nikodym derivative
  \[
  \frac{\textup{d}\nu}{\textup{d}\widetilde{\nu}\circ g_v^{-1}}\bigg{\rvert}_{X_v}
	=\lvert g_v'\circ g_v^{-1}\rvert^{\mdim}\bigg{\rvert}_{X_v}\cdot\bigg(\sum_{v'\in V}\int_{Y_{v'}} \lvert g_{v'}'\rvert^{\mdim}\textup{d}\widetilde{\nu}\bigg)^{-1}.
  \]
\item $L_{\om}^{v,j}=g_{i(\om_1)}(\widetilde{L}_{\om}^{v,j})$ for $v\in V$, $j\in\{1,\ldots,n_v\}$ and $\om\in T_v^*$. Define a function $f\colon E_A^{\infty}\to\mathbb R$ by $f(\om)\defeq\lvert g'_{i(\om_1)}\circ\widetilde{\bij}(\om)\rvert^{\mdim}$. Since $\lvert\widetilde{L}_{\om}^{v,j}\rvert=r_{\om}\lvert\widetilde{L}^{v,j}\rvert$, we have 
	\begin{align*}
  \lim_{n\to\infty}\sum_{v\in V}\sum_{j=1}^{n_v}\sum_{\om\in T_v^n}\lvert L_{\om}^{v,j}\rvert^{\mdim}
  &\ \, =\ \ \lim_{n\to\infty}\sum_{v\in V}\sum_{j=1}^{n_v}\sum_{\om\in T_v^n}\left(r_{\om}\lvert\widetilde{L}^{v,j}\rvert\cdot\lvert g'_{i(\om_1)}(x^{\om})\rvert\right)^{\mdim}\\
	&\ \, =\ \ \lim_{n\to\infty}\sum_{v\in V}\sum_{j=1}^{n_v}\lvert \widetilde{L}^{v,j}\rvert^{\mdim}\mathcal L_{-\mdim\ze}^n(f)(\om^v)\\
  &\stackrel{\eqref{eq:convergenceperron}}{=}\sum_{v\in V}\sum_{j=1}^{n_v}\lvert \widetilde{L}^{v,j}\rvert^{\mdim}\eigenf_{-\mdim\ze}(\om^v)\cdot\left(\sum_{v'\in V}\int_{Y_{v'}}\lvert g'_{v'}\rvert^{\mdim}\textup{d}\widetilde{\nu}\right),
  \end{align*}
  where $x^{\om}\in\widetilde{\bij}[\om]$ for each $\om\in E_A^*$ and $\om^v\in I_v^{\infty}$ for $v\in V$. Note that the above equation can be rigorously proven by using the Bounded Distortion Lemma (Lem.~\ref{lem:bd}).
\item\label{en:conformal} From the fact that $(R_e)_{e\in E}$ are contractions and each $g_v'$ is H\"older continuous and bounded away from zero, one can deduce that there exists a cGDS consisting of iterates of $\Phi\defeq(\phi_e)_{e\in E}$ which all are contractions. As this iterate also generates $F$, it follows that $F$ is a limit set of a cGDS.
\end{enumerate}

\begin{proof}[Proof of Thm.~\ref{thm:imafcm}]
Using \ref{en:diffbar} to \ref{en:conformal} an application of Thm.~\ref{thm:curvatureresult}\ref{curvatureresult:average} and \ref{curvatureresult:nonlattice} to $F$ and of Thm.~\ref{thm:similars} to $K$ proves Thm.~\ref{thm:imafcm}\ref{imafcm_average} and \ref{imafcm_nonlattice}.

The structure of the proof of \ref{imafcm_lattice} is taken from the proof of \cite[Thm.~2.17]{KesKom}. Here, the definitions of $\Delta$ and $(B(\kappa),f_{\kappa})$ are slightly different to the ones in \cite{KesKom}. We just provide the steps which require modification and refer the reader to \cite{KesKom} for detailed justifications. 

Write $R\eqdef(R_e)_{e\in E}$ and let $r_e\in(0,1)$ denote a Lipschitz constant of $R_e$ for $e\in E$. Note that $g_v\colon Y_v\to X_v$ is bijective for every $v\in V$. For $e\in E$ define 
\[
\phi_e\defeq g_{i(e)}\circ R_e\circ g_{t(e)}^{-1}
\]
and set $\Phi\defeq(\phi_e)_{e\in E}$. As is justified in \cite[p.\,2506]{KesKom} we can assume without loss of generality that $\phi_e$ are contractions. Then $\Phi$ is a cGDS and $F$ is the associated limit set.
The code space associated with $\Phi$ is also $E_A^{\infty}$. We let $\widetilde{\bij}$ and $\bij$ respectively denote the code maps from $E_A^{\infty}$ to $K$ and $F$. They satisfy $\bij(\om)=g_{i(\om_1)}\circ\widetilde{\bij}(\om)$ for $\om\in E_A^{\infty}$.
For a fixed $n\in\mathbb N\cup\{0\}$ define
\begin{align*}
  &\Delta_n\defeq\Big\{\bigcup_{i=1}^l[\kappa^{(i)}]\ \Big{\vert}\ \kappa^{(i)}\in E_A^{n},\, l\in\{1,\ldots,\card E_A^n\},\, \bigcup_{i=1}^l\langle\widetilde{\bij}[\kappa^{(i)}]\rangle\ \text{is an interval},\\
  &\hphantom{\Delta_n\defeq\Big\{\bigcup_{i=1}^l[\kappa^{(i)}]\mid\ }
  \bigcup_{i=1}^l\widetilde{\bij}[\kappa^{(i)}]\cap\widetilde{\bij}[\om]=\varnothing\ \text{for every}\ \om\in E_A^n\setminus\{\kappa^{(1)},\ldots,\kappa^{(l)}\}\Big\}.
\end{align*}
(Note that if the strong separation condition was satisfied, then $\Delta_n=\{[\om]\mid\om\in E_A^n\}$.)
Note that $\Delta_n\neq\varnothing$ for all $n\in\mathbb N$ because of the OSC and set $\Delta\defeq\bigcup_{n\in\mathbb N\cup\{0\}}\Delta_n$. 
Now, fix an $n\in\mathbb N\cup\{0\}$ and a $\kappa=\bigcup_{i=1}^l[\kappa^{(i)}]\in\Delta_n$ and choose $\theta>0$ such that $\bigcup_{i=1}^l\langle\widetilde{\bij}[\kappa^{(i)}]\rangle_{2\theta}\cap \widetilde{\bij}[\om]=\varnothing$ for every $\om\in E_A^n\setminus\{\kappa^{(1)},\ldots,\kappa^{(l)}\}$.
Then $B(\kappa)\defeq \bigcup_{i=1}^l\langle\widetilde{\bij}[\kappa^{(i)}]\rangle_{\theta}$ is a non-empty Borel subset of $\mathbb R$ satisfying $F_{\eps}\cap B(\kappa)=(F\cap B(\kappa))_{\eps}$ for all $\eps<\theta$. 

Let $\{L^{v,j}\}_{v\in V,\,j\in\{1,\ldots,n_v\}}$ denote the primary gaps of $F$ and $\{L_{\om}^{v,j}\}_{v\in V, j\in\{1,\ldots,n_v\}}$ the associated \main\ gaps.	
For constructing the function $f_{\kappa}$ fix an $m\in\mathbb N$ and choose $M\in\mathbb N$ so that 
\renewcommand{\theenumi}{(\roman{enumi})}
\begin{enumerate}
\item $\ee^{-M}<\theta$ and that 
\item $\lvert L_{\om}^{v,j}\rvert>2\ee^{-M}$ holds for all $v\in V$, $j\in\{1,\ldots,n_v\}$ and $\om\in T_v^m$ for which $L_{\om}^{v,j}\subset B(\kappa)$.
\end{enumerate}
Then for all $T\geq M$ we have
\begin{align}
  \leb^0\left(\partial F_{\ee^{-T}}\cap B(\kappa)\right)/2
  &= \sum_{v\in V}\sum_{j=1}^{n_v}\card\{\om\in T_v^*\mid L_{\om}^{v,j}\subseteq B(\kappa),\ \lvert L_{\om}^{v,j}\rvert>2\ee^{-T}\}+1\nonumber\\
  &\leq \sum_{v\in V}\sum_{j=1}^{n_v}\sum_{\om\in T_v^{m}}\Xi_{\om}^{v,j}(\ee^{-T})+\underbrace{\sum_{v\in V}n_v\sum_{j=1}^{m-1}(\card E)^{j-1}+1}_{\mbox{\normalsize $\eqdef c_m$}},\label{eq:nonexupper}
\end{align}
where 
\begin{equation*}
  \Xi_{\om}^{v,j}(\ee^{-T})\defeq\card\{\omneu\in T_{i(\om_1)}^*\mid L_{\omneu\om}^{v,j}\subseteq B(\kappa),\ \lvert L_{\omneu\om}^{v,j}\rvert> 2\ee^{-T}\}.
\end{equation*}
Likewise
\begin{equation*}
  \leb^0\left(\partial F_{\ee^{-T}}\cap B(\kappa)\right)/2
  \geq \sum_{v\in V}\sum_{j=1}^{n_v}\sum_{\om \in T_v^m}\Xi_{\om}^{v,j}(\ee^{-T}).
\end{equation*}
For providing bounds on $\Xi_{\om}^{v,j}(\ee^{-T})$, we let $\xi$ and $\ze$ denote the geometric potential functions associated with $\Phi$ and $R$. 
For $\om\in E_A^{\infty}$ we have
\begin{align*}
  \xi(\om)
  &=\hspace{-0.05cm}-\ln\lvert\phi'_{\om_1}(\bij\sigma \om)\rvert\\
  &=\hspace{-0.05cm}-\ln \lvert g'_{i(\om_1)}(R_{\om_1}g_{t(\om_1)}^{-1}\bij\sigma \om)\rvert\hspace{-0.05cm} -\hspace{-0.05cm}\ln\lvert R'_{\om_1}(g_{t(\om_1)}^{-1}\bij\sigma \om)\rvert\hspace{-0.05cm} +\hspace{-0.05cm}\ln \lvert g_{t(\om_1)}'(g_{t(\om_1)}^{-1}\bij\sigma \om)\rvert\\
  &=\hspace{-0.05cm}-\ln \lvert g_{i(\om_1)}'(\widetilde{\bij} \om)\rvert+\ze(\om)+\ln \lvert g'_{t(\om_1)}(\widetilde{\bij}\sigma \om)\rvert.
\end{align*}
Therefore, $\psi\colon E_A^{\infty}\to\mathbb R$ given by $\psi(\om)\defeq-\ln \lvert g'_{i(\om_1)}(\widetilde{\bij} \om)\rvert$ defines a function lying in $\mathcal C(E_A^{\infty})$ which satisfies 
\[
\xi-\ze=\psi-\psi\circ\sigma.
\]
Let $c$ be the common H\"older constant of $g_v$ for $v\in V$ and let $k>0$ be such that for each $v\in V$ we have that $\lvert g'_v\rvert\geq k$ on $W_v$. Then for all $x,y\in\langle \widetilde{\bij}[\om]\rangle$, where $\om\in I_v^n$ for $n\in\mathbb N$ and $v\in V$ we have that 
\begin{equation}\label{eq:bdg}
  \left\lvert\frac{g'_v(x)}{g'_v(y)}\right\rvert
  \leq \left\lvert\frac{g'_v(x)-g'_v(y)}{g'_v(y)}\right\rvert+1
  \leq \frac{c\lvert x-y\rvert^{\alpha}}{k}+1
  \leq \max_{\om\in I_v^n}\frac{c\lvert\langle \bij[\om]\rangle\rvert^{\alpha}}{k}+1
	\eqdef \bdneu_n.
\end{equation} 
Clearly, $\bdneu_n\to 1$ as $n\to\infty$. We let $\om^v\in I_v^{\infty}$ be arbitrary and $\om\in T_v^m$. Then 
\begin{align*}
  \lvert L_{\omneu\om}^{v,j}\rvert
  &=\lvert g_{i(\omneu_1)}\widetilde{L}_{\omneu\om}^{v,j}\rvert
  \leq \lvert g'_{i(\omneu_1)}(R_{\omneu\om}\widetilde{\bij}\om^v)\rvert\bdneu_m\cdot\lvert R'_{\omneu}(R_{\om}\widetilde{\bij}\om^v)\rvert\cdot\lvert\widetilde{L}_{\om}^{v,j}\rvert\\
  &=\lvert (g_{i(\omneu_1)}\circ R_{\omneu})'(R_{\om}\widetilde{\bij}\om^v)\rvert\cdot\bdneu_m\lvert \widetilde{L}_{\om}^{v,j}\rvert\\
  &= \lvert\phi'_{\omneu}(\phi_{\om}\bij\om^v)\rvert\cdot\lvert g'_{i(\om_1)}(R_{\om}\widetilde{\bij}\om^v)\rvert\cdot\bdneu_m\cdot r_{\om}\lvert \widetilde{L}^{v,j}\rvert\\
  &=\exp\big(-S_{n(\omneu)}\xi(\omneu\om \om^v)-\psi(\om \om^v)+\ln(\bdneu_m\cdot r_{\om}\lvert \widetilde{L}^{v,j}\rvert)\big).
\end{align*}
Therefore, for such $\om^v\in I_v^{\infty}$, $T\geq\max\{M,\widetilde{M}\}$ and $\om\in T_v^m$ we have that
\begin{align*}
&\Xi_{\om}^{v,j}(\ee^{-T})
\leq\card\{\omneu\in T_{i(\om_1)}^*\mid L_{\omneu\om}^{v,j}\subseteq B(\kappa),\\
&\hphantom{\Xi_{\om}^{v,j}(\ee^{-T})\leq\card\{\omneu\in T_{i(\om_1)}^*\mid\ } S_{n(\omneu)}\xi(\omneu\om \om^v)< T+\ln(\frac{\bdneu_{m} r_{\om}\lvert\widetilde{L}^{v,j}\rvert}{2})-\psi(\om \om^v)\}.
\end{align*}
Applying \cite[Thm.~3.9 and Rem.~3.10]{KesKom} yields
\begin{align}
  &\leb^0(\partial F_{\ee^{-T}}\cap B(\kappa))/2-c_m\nonumber\\
  &\leq \sum_{v\in V}\sum_{j=1}^{n_v}\sum_{\om\in T_v^m}\sum_{n=0}^{\infty}\sum_{\omneu:\sigma^n \omneu=\om \om^v}\mathds 1_{\kappa}(\omneu)\cdot\mathds 1_{\big{\{}S_{n}\xi(\omneu)< T+\ln(\bdneu_{m}r_{\om}\lvert\widetilde{L}^{v,j}\rvert/2)-\psi(\om \om^v)\big{\}}}\nonumber\\
  &\sim \sum_{v\in V}\sum_{j=1}^{n_v}\sum_{\om\in T_v^m}\hspace{-0.1cm}\frac{\aaa\eigenf_{-\mdim\ze}(\om \om^v)\int_{\kappa}\hspace{-0.05cm}\ee^{-\mdim\aaa\left\lceil\hspace{-0.05cm}\frac{\psi(\omneu)-\psi(\om \om^v)}{\aaa}+\frac{1}{\aaa}\ln{\frac{2\ee^{-T}}{\bdneu_{m}r_{\om}\lvert\widetilde{L}^{v,j}\rvert}}+\frac{\psi(\om \om^v)}{\aaa}\hspace{-0.05cm}\right\rceil}\hspace{-0.05cm}\textup{d}\nu_{-\mdim\ze}(\omneu)}{\left(1-\ee^{-\mdim\aaa}\right)\int\ze\textup{d}\mu_{-\mdim\ze}}	\label{long}
\end{align}
  Define $U\defeq\aaa\left(1-\ee^{-\mdim\aaa}\right)^{-1}\left(\int\ze\textup{d}\mu_{-\mdim\ze}\right)^{-1}$. Using that $\ln r_{\om}\in\aaa\mathbb Z$ for every $\om\in E_A^{*}$, the right hand side of \eqref{long} can be rewritten as 
  \begin{equation*}
    \sum_{v\in V}\sum_{j=1}^{n_v}\sum_{\om\in T_v^m}Ur_{\om}^{\mdim}\eigenf_{-\mdim\ze}(\om \om^v)\int_{\kappa} \ee^{-\mdim\aaa\left\lceil\frac{\psi(\omneu)}{\aaa}+\frac{1}{\aaa}\ln{\frac{2\ee^{-T}}{\bdneu_{m}\lvert\widetilde{L}^{v,j}\rvert}}\right\rceil}\textup{d}\nu_{-\mdim\ze}(\omneu).
	\end{equation*}
  Defining the function $f_{\kappa}\colon\mathbb R^+\to\mathbb R^+$ by
  \[
  f_{\kappa}(T)\defeq \ee^{-\mdim T}\sum_{v\in V}\sum_{j=1}^{n_v}\sum_{\om\in T_v^m} U \eigenf_{-\mdim\ze}(\om \om^v)r_{\om}^{\mdim}\int_{\kappa} \ee^{-\mdim\aaa\left\lceil\frac{\psi(\omneu)}{\aaa}+\frac{1}{\aaa}\ln{\frac{2\ee^{-T}}{\lvert\widetilde{L}^{v,j}\rvert}}\right\rceil}\textup{d}\nu_{-\mdim\ze}(\omneu)
  \]
  we thus have that for all $t>0$ there exists a $\widetilde{M}\in\mathbb R$ such that
  \begin{equation*}
    \ee^{-\mdim T}\leb^0(\partial F_{\ee^{-T}}\cap B(\kappa))/2
    \leq (1+t)\bdneu_m^{\mdim}f_{\kappa}(T+\ln\bdneu_m)+c_m\ee^{-\mdim T}.
  \end{equation*}
  for all $T\geq\widetilde{M}$ and likewise, 	
  \begin{equation*}
    \ee^{-\mdim T}\leb^0(\partial F_{\ee^{-T}}\cap B(\kappa))/2
    \geq (1-t)\bdneu_m^{-\mdim}f_{\kappa}(T-\ln\bdneu_m).
  \end{equation*}
  Clearly, $f_{\kappa}$ is periodic with period $\aaa$. Thus, \cite[Lem.~5.1 (ii)]{KesKom} is satisfied for $B=B(\kappa)$ and $f=f_{\kappa}$.
  For showing validity of \cite[Lem.~5.1 (i)]{KesKom}, we set $\underline{\beta}\defeq \min\{\{\aaa^{-1}\ln\lvert\widetilde{L}^{v,j}\rvert\}\mid v\in V, j\in\{1,\ldots,n_v\}\}$ and $\overline{\beta}\defeq\max\{\{\aaa^{-1}\ln\lvert\widetilde{L}^{v,j}\rvert\}\mid v\in V, j\in \{1,\ldots,n_v\}\}$. We first assume that $\underline{\beta}>0$ and consider the following four cases, where we let $q^*\in\mathbb N\cup\{0\}$ be maximal such that $\underline{\beta}+q^*(1-\overline{\beta})\leq\overline{\beta}$.

  \textsc{Case 1}: $\displaystyle\underline{D}\defeq\{\om\in E_A^{\infty}\mid\{\aaa^{-1}\psi(\om)\}<\underline{\beta}\}\not=\varnothing$.
  
  \textsc{Case 2}: $\displaystyle\overline{D}\defeq\{\om\in E_A^{\infty}\mid\{\aaa^{-1}\psi(\om)\}>\overline{\beta}\}\not=\varnothing$.  
  
  \textsc{Case 3}:
  There exists a $q\in\{0,\ldots,q^*\}$ such that\\
  \hphantom{\textsc{Case 3}:\ }
  \indent $\displaystyle D_q\defeq \{\om\in E_A^{\infty}\mid\underline{\beta}+q(1-\overline{\beta})<\{\aaa^{-1}\psi(\om)\}<\underline{\beta}+(q+1)(1-\overline{\beta})\}\neq\varnothing$.  
   
  \textsc{Case 4}:
  $\displaystyle\{\om\in E_A^{\infty}\mid\{\aaa^{-1}\psi(\om)\}\subseteq\{\underline{\beta}+q(1-\overline{\beta})\mid q\in\{0,\ldots,q^*\}\}\}= E_A^{\infty}$.

 Note that Case 4 obtains if neither of the cases 1-3 obtains.
 With the same methods as in \cite[p.\,2508 f.]{KesKom}, in particular using the same functions $T_n$, one can deduce that $f_{\kappa}$ is not equal to an almost everywhere constant function in all four cases.
  The conclusion of the proof is the same as in \cite[proof of Thm.~2.17]{KesKom}.
 \end{proof}

\begin{proof}[Proof of Thm.~\ref{thm:analytic}]
	We define an operator $\widetilde{\mathcal L}\colon \mathcal C(\mathbb R)\to\mathcal C(\mathbb R)$ by setting 
	\[
		\widetilde{\mathcal L}(g)(x)
		\defeq \sum_{e\in T_v}\lvert\phi'_e(x)\lvert^{\mdim}\cdot g\circ\phi_e(x)
	\]
	for $x\in X_v$ and $v\in V$.
	Letting $\xi$ denote the geometric potential function associated with $\Phi$ and letting $\bij$ denote the code map from the code space $E_A^{\infty}$ to $F$, we see that $\widetilde{\mathcal L}(g)(\bij\om)=\mathcal L_{-\mdim\xi}(g\circ\bij)(\om)$, where $\mdim$ denotes the Minkowski dimension of $F$.
	
	Since the maps $\phi_e$ are analytic, there exist open neighbourhoods $W_v\supset X_v$ of $X_v$ in $\mathbb C$ on which the maps $\phi_e$ are analytic for $e\in T_v$. By \cite[Lem.~4.2.12]{Urbanski_Buch} the functions $\widetilde{\mathcal L}^n \mathbf{1}\rvert_{W_v}$ are uniformly bounded and the bound is independent of $n\in\mathbb N$. Thus, for $v\in V$, $\widetilde{\mathcal L}^n \mathbf{1}\colon W_v\to\mathbb C$ form a normal family in the sense of Montel. Here, $\mathbf{1}$ denotes the constant one-function on $E_A^{\infty}$.
	By \eqref{eq:convergenceperron} we have that $\widetilde{\mathcal L}^n \mathbf{1}\circ\bij$ converges uniformly to $\eigenf_{-\mdim\xi}$ on $E_A^{\infty}$. Therefore, $\widetilde{\mathcal L}^n\mathbf{1}\rvert_{W_v}$ converges to an analytic extension of $\eigenf_{-\mdim\xi}$ on $W_v$. We denote this analytic extension by $\eigenf^v$ and set $\widetilde{\psi}_v\defeq\mdim^{-1}\ln \eigenf^v$.	
	Since $\xi$ is lattice, there exist $\ze,\psi\in\mathcal C(E_A^{\infty})$ such that
	\[
		\xi-\ze=\psi-\psi\circ\sigma
	\]
	and such that the range of $\ze$ is contained in a discrete subgroup of $\mathbb R$. We let $\aaa>0$ denote the maximal real number such that $\ze(E_A^{\infty})\subset\aaa\mathbb Z$. Note that $\widetilde{\psi}_v$ satisfies
	\[
	\widetilde{\psi}_v\circ\bij\rvert_{I_v}
	=\psi\rvert_{I_v}+\mdim^{-1}\ln\eigenf_{-\mdim\ze}\rvert_{I_v},
	\]
	as $\eigenf^v$ satisfies $\eigenf^v\circ\bij\rvert_{I_v}=\ee^{\mdim\psi}\eigenf_{-\mdim\ze}\rvert_{I_v}$, where we used that 
	\[
	\mathcal{L}_{-\mdim\xi}(\ee^{\mdim\psi}\eigenf_{-\mdim\ze})(x)
	=\sum_{y:\sigma y=x}\ee^{\mdim(\psi-\xi)(y)}\eigenf_{-\mdim\ze}(y)
	=\ee^{\mdim\psi(x)}\mathcal{L}_{-\mdim\ze}(\eigenf_{-\mdim\ze})(x)
	=\ee^{\mdim\psi}\eigenf_{-\mdim\ze}(x)
	\]
	which implies $\eigenf_{-\mdim\xi}=\ee^{\mdim\psi}\eigenf_{-\mdim\ze}$.
	We define $X_v\eqdef[a_v,b_v]$ for $v\in V$ and introduce the functions $\widetilde{g}_v\colon[a_v,b_v]\to\mathbb R$ given by
	\[
	\widetilde{g}_v(x)\defeq\int_{a_v}^x \ee^{\widetilde{\psi}_v(y)}\textup{d} y/D_v+2v,
	\]
	where $D_v\defeq\int_{a_v}^{b_v} \ee^{\widetilde{\psi}_v(y)}\textup{d} y$. Notice $\widetilde{g}_v([a_v,b_v])=[2v,2v+1]$.
	As $\widetilde{\psi}^v$ is analytic by definition, the fundamental theorem of calculus implies that $\widetilde{g}_v'(x)=\ee^{\widetilde{\psi}_v(x)}/D_v$, giving
	\[
	\ln\widetilde{g}'_v=\widetilde{\psi}_v-\ln D_v.
	\]
	Furthermore, the analyticity of $\widetilde{\psi}_v$ implies that $\widetilde{\psi}_v$ is bounded on $X_v$. Therefore, $\widetilde{g}'_v$ is bounded away from both 0 and $\infty$ and hence $\widetilde{g}_v$ is invertible. Set
	\[
	g_v\colon[2v,2v+1]\to[a_v,b_v],\quad g_v\defeq\widetilde{g}_v^{-1}
	\]
	and extend $g_v$ to an analytic function on an open neighbourhood $\mathcal U_v$ of the interval $[2v,2v+1]$ such that $\lvert g'_v\rvert>0$ on $\mathcal U_v$. For $e\in E$ we define 
	\[
	R_e\defeq g_{i(e)}^{-1}\circ\phi_e\circ g_{t(e)}
	\]
	and introduce the code map $\widetilde{\bij}$ given by $\widetilde{\bij}\rvert_{I_v}\defeq g_v^{-1}\circ\bij$ for $v\in V$. For $\om\in E_A^{\infty}$ we then have
	\begin{align*}
	  &-\ln R'_{\om_1}(\widetilde{\bij}\sigma\om)\\
	  &\qquad=-\ln \widetilde{g}'_{i(\om_1)}(\phi_{\om_1}g_{t(\om_1)}\widetilde{\bij}\sigma\om)-\ln\phi'_{\om_1}(g_{t(\om_1)}\widetilde{\bij}\sigma\om)+\ln\widetilde{g}'_{t(\om_1)}(g_{t(\om_1)}\widetilde{\bij}\sigma\om)\\
	  &\qquad=-\widetilde{\psi}_{i(\om_1)}(\bij\om)+\xi(\om)+\widetilde{\psi}_{t(\om_1)}(\bij\sigma\om)+\ln D_{i(\om_1)}-\ln D_{t(\om_1)}\\
	  &\qquad=-\psi(\om)-\mdim^{-1}\ln(\eigenf_{-\mdim\ze}(\om)/\eigenf_{-\mdim\ze}(\sigma\om))+\psi(\sigma\om)+\xi(\om)+\ln(D_{i(\om_1)}/D_{t(\om_1)})\\
	  &\qquad=\ze(\om)-\mdim^{-1}\ln\frac{\eigenf_{-\mdim\ze}(\om)}{\eigenf_{-\mdim\ze}(\sigma\om)}+\ln\frac{D_{i(\om_1)}}{D_{t(\om_1)}}.
	\end{align*}
	Since the range of $\ze$ is contained in the group $\aaa\mathbb Z$ and $\xi$ and $\psi$ are bounded on $E_A^{\infty}$, $\ze$ in fact takes a finite number of values. The continuity of $\ze$ implies that there exists an $M\in\mathbb N$ such that $\ze$ is constant on each cylinder set $[\om]$ for $\om\in E_A^M$. This clearly implies that $\mathcal L_{-\mdim\ze}^n \mathbf{1}$ is constant on $[\om]$ for all $\om\in E_A^M$ and all $n\in\mathbb N$. Thus, \eqref{eq:convergenceperron} implies that also $\eigenf_{-\mdim\ze}$ is constant on cylinder sets of length $M$. This can be seen by considering $\lvert\eigenf_{-\mdim\ze}(\om)-\eigenf_{-\mdim\ze}(\omneu)\rvert$ for $\omneu,\om$ lying in the same cylinder set of length $M$ and applying the triangle inequality. Therefore, $\om\mapsto-\ln\lvert R'_{\om_1}(\widetilde{\bij}\sigma\om)\rvert$ is constant on cylinder sets of length $M+1$. 
	Since for each $\om\in E_A^{M+1}$ the set $\{\widetilde{\bij}\omneu\mid\omneu\in[\om]\}$ has accumulation points and is compact and the map $R'_e$ is analytic by construction, it follows that $R'_e$ is constant on its domain of definition. Therefore, the maps $R_e$ are similarities. From the fact that $\phi_e$ are contractions and each of the $g_v'$ is differentiable and bounded away from zero, one can deduce that there exists an iterate $\widetilde{R}$ of $R\defeq(R_e)_{e\in E}$ which solely consists of contractions and thus is an sGDS.  
 \end{proof}


\begin{thebibliography}{7}

\bibitem{Beardon}
A.~F.~Beardon, \emph{The geometry of discrete groups}, Graduate Texts in
  Mathematics, vol.~91, Springer, New York 1995.
\newblock Corrected reprint of the 1983 original.

\bibitem{Bohl}
T.~J.~Bohl, Fractal curvatures and Minkowski content of self-conformal sets, preprint arXiv:1211.3421.

\bibitem{Bowen_equilibrium}
R.~Bowen, \emph{Equilibrium states and the ergodic theory of {A}nosov
  diffeomorphisms}, Lecture Notes in Mathematics, vol.~470, revised edn., Springer, Berlin 2008.

\bibitem{Connes}
A.~Connes, \emph{Noncommutative geometry}, Academic Press Inc., San Diego, CA, 1994.

\bibitem{Trken_graph-directed}
B.~Demir, A.~Deniz, \c{S}.~Ko\c{c}ak, A.~E.~\"Ureyen, Tube formulas for graph-directed fractals, \emph{Fractals} \textbf{18} (2010), 349--361.

\bibitem{Trken}
A.~Deniz, \c{S}.~Ko\c{c}ak, Y.~\"Ozdemir, A.~V.~Ratiu, A.~E.~\"Ureyen, On the
  {M}inkowski measurability of self-similar fractals in {$R^d$}, \emph{Turkish Journal of Mathematics} \textbf{37} (2013), 830--846.

\bibitem{DenkerUrbanski}
M.~Denker, M.~Urba{\'n}ski, On the existence of conformal measures, \emph{Trans. Amer. Math. Soc.} (2) \textbf{328} (1991), 563--587.

\bibitem{Falconer_Minkowski}
K.~J.~Falconer, On the {M}inkowski measurability of fractals.
\emph{Proc. Amer. Math. Soc.} (4) \textbf{123} (1995), 1115--1124.

\bibitem{Falconer_Foundation}
K.~J.~Falconer, \emph{Fractal geometry. {M}athematical foundations and applications},
  second edn., John Wiley \& Sons Inc., Hoboken, NJ, 2003.

\bibitem{Falconer_Samuel}
K.~J.~Falconer, T.~Samuel, Dixmier traces and coarse multifractal analysis,
\emph{Ergodic Theory Dynam. Systems} (2) \textbf{31} (2011), 369--381.

\bibitem{Uta}
U.~Freiberg, S.~Kombrink, Minkowski content and local {M}inkowski content for
  a class of self-conformal sets,
\emph{Geom. Dedicata} \textbf{159} (2012), 307--325.

\bibitem{Gatzouras}
D.~Gatzouras, Lacunarity of self-similar and stochastically self-similar sets,
\emph{Trans. Amer. Math. Soc.} (5) \textbf{352} (2000), 1953--1983.

\bibitem{Guido_Isola}
D.~Guido, T.~Isola, Dimensions and singular traces for spectral triples, with
  applications to fractals,
\emph{J. Funct. Anal.} (2) \textbf{203} (2013), 362--400.

\bibitem{KesKom}
M.~Kesseb{\"o}hmer, S.~Kombrink, Fractal curvature measures and {M}inkowski
  content for self-conformal subsets of the real line,
 \emph{Adv. Math.} (4-6) \textbf{230} (2012), 2474--2512.

\bibitem{survey}
S.~Kombrink, {A} {S}urvey on {M}inkowski measurability of {S}elf-{S}imilar and {S}elf-{C}onformal {F}ractals in {$R^d$},
in: \emph{{F}ractals in {P}ure {M}athematics}, {C}ontemporary {M}athematics 600, {A}mer. {M}ath. {S}oc, 2013, 135--159.

\bibitem{Diss}
S.~Kombrink, Fractal curvature measures and {M}inkowski content for limit sets of conformal function systems. Dissertation, Bremen, Germany, 2011,
  http://nbn-resolving.de/urn:nbn:de:gbv:46-00102477-19.

\bibitem{Lalley}
S.~P.~Lalley, Renewal theorems in symbolic dynamics, with applications to
  geodesic flows, non-{E}uclidean tessellations and their fractal limits,
\emph{Acta Math.} (1-2) \textbf{163} (1989), 1--55.

\bibitem{Dundee}
M.~L.~Lapidus, Vibrations of fractal drums, the {R}iemann hypothesis, waves in fractal media and the {W}eyl-{B}erry conjecture,
in: \emph{Ordinary and partial differential equations}, {V}ol.~{IV}
  ({D}undee, 1992), Pitman Res. Notes Math. Ser., vol.~289, pp.
  126--209. Longman Sci. Tech., Harlow (1993).

\bibitem{Lapidus_Frankenhuysen_Springer}
M.~L.~Lapidus, M.~van Frankenhuijsen, \emph{Fractal geometry, complex dimensions and zeta functions. {G}eometry and spectra of fractal strings},
Springer Monographs in Mathematics, Springer, New York, 2006.

\bibitem{LapPeaWin}
M.~L.~Lapidus, E.~Pearse, S.~Winter, Minkowski measurability results for
  self-similar tilings and fractals with monophase generators,
in: \emph{{F}ractals in {P}ure {M}athematics}, {C}ontemporary {M}athematics 600, {A}mer. {M}ath. {S}oc, 2013, 185--203.

\bibitem{LapPom}
M.~L.~Lapidus, C.~Pomerance, The {R}iemann zeta-function and the
  one-dimensional {W}eyl-{B}erry conjecture for fractal drums,
\emph{Proc. London Math. Soc.} (3) \textbf{66}(1) (1993), 41--69.

\bibitem{WinterLlorente}
M.~Llorente, S.~Winter, A notion of {E}uler characteristic for fractals,
\emph{Math. Nachr.} (1-2) \textbf{280} (2007), 152--170.

\bibitem{Urbanski_Buch}
R.~D.~Mauldin, M.~Urba{\'n}ski, \emph{Graph directed {M}arkov systems. {G}eometry and dynamics of limit sets}, Cambridge Tracts in Mathematics, vol.~148,
Cambridge University Press, Cambridge, 2003.

\bibitem{Nicholls}
P.~J.~Nicholls, \emph{The ergodic theory of discrete groups}, London
  Mathematical Society Lecture Note Series, vol.~143,
Cambridge University Press, Cambridge, 1989.

\bibitem{Patterson}
S.~J.~Patterson, The limit set of a {F}uchsian group,
\emph{Acta Math.} (3-4) \textbf{136} (1976), 241--273.

\bibitem{PokornyWinter}
D.~Pokorn{\'y}, S.~Winter, Scaling exponents of curvature measures, \emph{J. Fractal Geom.} \textbf{1(2)} (2014), 177--219.

\bibitem{RatajWinter}
J.~Rataj, S.~Winter, On volume and surface area of parallel sets,
\emph{Indiana Univ. Math. J.} (5) \textbf{59} (2010), 1661--1685.

\bibitem{Ruelle_gas}
D.~Ruelle, Statistical mechanics of a one-dimensional lattice gas,
\emph{Comm. Math. Phys.} \textbf{9} (1968), 267--278.

\bibitem{Samuel}
T.~Samuel,\,A commutative noncommutative fractal geometry.\,Ph.D.\,Thesis, {S}t\,{A}ndrews,\,2010.

\bibitem{Sullivan}
D.~Sullivan, The density at infinity of a discrete group of hyperbolic
  motions, \emph{Inst. Hautes \'Etudes Sci. Publ. Math.} (50) (1979), 171--202.

\bibitem{Walters_convergence}
P.~Walters, Convergence of the {R}uelle operator for a function satisfying
  {B}owen's condition,
\emph{Trans. Amer. Math. Soc.} (1) \textbf{353} (2001), 327--347.

\bibitem{Winter_thesis}
S.~Winter, Curvature measures and fractals, 
\emph{Dissertationes Math. (Rozprawy Mat.)} \textbf{453} (2008).



\end{thebibliography}
\end{document}